\documentclass[11pt]{article}
\small
\normalsize
\usepackage[latin1]{inputenc}
\usepackage[T1]{fontenc}
\usepackage{amsmath}
\usepackage{amssymb}
\usepackage{amsfonts}
\usepackage{amsbsy}
\usepackage{amscd}
\usepackage{theorem}
\usepackage{epsf}
\usepackage{psfrag}
\allowdisplaybreaks

\def\build#1_#2^#3{\mathrel{\mathop{\kern 0pt#1}\limits_{#2}^{#3}}}
\def\noi{{\noindent}}
\def\be{\begin{equation}}
\def\ee{\end{equation}}
\def\ba{\begin{eqnarray*}}
\def\ea{\end{eqnarray*}}

\def\E{{\bf E}}
\def\P{\hbox{\bf P}}
\def\cqfd{ \hfill $\blacksquare$ }

\def\o{\omega}
\def\OO{\Omega }

\def\F{\hbox{\bf F}}

\def\f{{\cal F}}

\def\i{{\cal I}}

\def\n{{\cal N}}

\def\u{{\cal U}}
\def\U{{\cal U}}
\def\PP{{\bf P}}
\def\EE{{\bf E}}

\def\wh{\widehat}

\def\la{\longrightarrow}

\def\eps{\epsilon}
\def\ddd{\overline{d}}

\def\gg{\underline{g}}
\def\ggg{\overline{g}}
\def\locc{L^{*-1}}
\def\loc{L^{-1}}
\def\pw{\cdot \wedge }
\def\Yd{\underrightarrow{Y}}
\def\Yg{\underleftarrow{Y}}
\def\spos{\sigma^{\uparrow}}
\def\tpos{\tau^{\uparrow}}
\def\ipos{\cal{I}^{\uparrow}}

\def\tneg{\tau^{\downarrow}}
\def\Xupd{\overrightarrow{X}}
\def\Xupg{\overleftarrow{X}}
\def\Xd{\underrightarrow{X}}
\def\XdT{\underrightarrow{X}^{\scriptscriptstyle{T}} }
\def\XgT{\underleftarrow{X}^{\scriptscriptstyle{T}} }
\def\Xg{\underleftarrow{X}}
\def\Xneg{X^{\downarrow }}
\def\Xpos{X^{\uparrow}}
\def\tpssig{\sigma_x (t) }
\def\xlast{\wh X^{\sigma_x (t) }}
\def\xsigma{ X_{\sigma_x (t) }}

\def\Zg{\underleftarrow{Z}}
\def\Zd{\underrightarrow{Z}}
\def\sd{\underrightarrow{\sigma }}
\def\td{\underrightarrow{\tau }}

\def\tdT{\underrightarrow{\tau }^{\scriptscriptstyle{T} }_x}
\def\L{\Lambda }
\def\l{\lambda }
\def\gi{\underline{g}}
\def\gs{\overline{g}}

\def\oo{\widehat{\omega}}
\def\mesu{{\rm m}}

\def\noi{\noindent}
\def\dem{\noindent{\bf Proof.} }
\def\dis{\displaystyle }

\def\rem{\noindent{\bf Remark.} }
\def\rems{\noindent{\bf Remarks.} }

\def\conss{\noindent{\bf Consequences.} }

\def\build#1_#2^#3{\mathrel{\mathop{\kern 0pt#1}\limits_{#2}^{#3}}}
\newcommand{\fr}[2]{\frac{\dis #1 }{\dis #2 }}
\newtheorem{theorem}{Theorem}[section]
\newtheorem{lemma}[theorem]{Lemma}
\newtheorem{proposition}[theorem]{Proposition}
\newtheorem{corollary}[theorem]{Corollary}

{\theorembodyfont{\rmfamily}}
{\theorembodyfont{\rmfamily}}

\newcommand{\RR}{\mathbb{R}}

\renewcommand{\P}{\mathbb{P}}

\newcommand{\un}{\boldsymbol{1}}
\begin{document}

\title{ \bf PATH DECOMPOSITIONS FOR \\
REAL LEVY PROCESSES}
\author{ by  {\sc  Thomas Duquesne,}\\
{\it Ecole Normale Sup\'erieure de Cachan, C.M.L.A. , } \\
{\it 61, avenue du President Wilson } \\
{\it 94230 CACHAN.} \\
e-mail: duquesne@cmla.ens-cachan.fr \\
phone: 01 47 40 27 88, fax: 01 47 40 59 01.}
\vspace{4mm}
\maketitle

\begin{abstract}

Let $X$ be a real L\'evy process and let $\Xpos $ be the process conditioned to stay positive. We assume that $ 0  $ is regular for $(-\infty , 0 )$
and $(0, +\infty ) $ with 
respect to $X$. Using elementary excursion theory arguments, we provide a simple probabilistic description of the reversed 
paths of $X$ and $\Xpos $ at their first hitting time of $ (x, +\infty )$  and last passage time of $ (-\infty , x ] $, 
on a fixed time interval $[0, t]$, for a positive level $x$. From these reversion formulas, we derive an extension to general L\'evy processes
of Williams' decomposition theorems, Bismut's decomposition of the excursion
above the infimum and 
also several relations involving the reversed excursion under the maximum.
\end{abstract}

\renewcommand{\abstractname}{R\'esum\'e}

\begin{abstract}
  Soit $X$ un processus de L\'evy et $\Xpos $ le m\^eme
processus conditionn\'e \`a rester positif. On suppose que $0$ est
r\'egulier pour $(-\infty , 0 )$ et $(0, +\infty ) $ par rapport \`a $X$. Par des arguments simples de th\'eorie des excursions, nous 
d\'ecomposons la loi des trajectoires de $X$ et $\Xpos $ retourn\'ees aux temps d'entr\'ee de $ (x, +\infty )$ et de sortie de 
$ (-\infty , x ] $. De ces formules de reversion, on d\'eduit une extension au cas des processus de L\'evy  g\'en\'eraux, des 
th\'eor\`emes de d\'ecomposition de Williams, du th\'eor\`eme de d\'ecomposition de Bismut de l'excursion au dessus du minimum, ainsi que
plusieurs relations faisant intervenir l'excursion sous le maximum retourn\'ee.
\end{abstract}

\section{Introduction}
Let $ (X_t)_{t\geq 0 } $ be a real L\'evy process, that is a real valued process with homogeneous and independent increments. The supremum 
(resp. infimum) of $X$ on the 
time interval $[0, t]$ is denoted by $ S_t $ (resp. $I_t $). 
We assume that $0$ is regular for $(0, +\infty )$ and 
$(-\infty , 0)$ with respect to $X$. A classical result says 
that $X-I$ (resp. $X-S$) is a strong Markov process for 
which $0$ is regular (see Bingham \cite{Bi1} or Bertoin \cite{Be} chapter 6 for a proof). 
Let us denote by $L$ (resp. $L^* $) 
the local time at $0$ of $X-I$ (resp. $X-S$): they
are uniquely defined up to a multiplicative constant and 
their normalization is specified in Proposition 
\ref{bismut}.

  As Rogers noticed in \cite{Ro84}, 
as soon as the L\'evy measure charges the positive numbers,
$X-S$ may hit zero by $X$ jumping across the level of its previous maximum. 
The classical It\^o excursion measure of $S-X$ loses the information about this jump. 
Let us introduce the relevant 
definition of the excursion measure under the maximum (resp. above the infimum) denoted by  $N^* $ (resp. $N$). It  has the property to
record the final jump of the excursion, which 
represents the amount the excursion overshoots when $X$ attains a new maximum (resp. infimum). Let 
$(g_i , d_i ), \; i\in \i$ 
(resp. $(g_j , d_j ) , \; j \in {\i}^* $) the excursion intervals of $X-I$ 
(resp. $S-X $) above $0$.
We define the excursions above the infimum
and under the supremum by 
$$\left\{ 
\begin{array}{ll}
\o^i (s) = X_{ (g_i + s)\wedge d_i } - X_{g_i } , \quad i\in \i\\
\o^j (s) =  X_{ (g_j + s)\wedge d_j } - X_{g_j } , \quad j\in \i^*
\end{array} \right. $$
\noi Then, the point measures 
$$ \sum_{i\in \i} \delta_{( L_{g_i } , \o^i )} \quad {\rm and } \quad  
\sum_{j\in \i^*} \delta_{( L^*_{g_j } , \o^j )} $$
\noi are distributed respectively as $\un_{ \{ l\leq \eta \} } \n (dld\o ) $ and 
 $\un_{ \{ l\leq \eta^* \} } \n^* (dld\o ) $, where $\n $ and $\n^* $ are 
Poisson measures with respective intensities $dlN(d\o )$ and   
$dlN^*(d\o )$, and where
$$\left\{
\begin{array}{ll}
\eta =\inf \left\{ t\geq 0 \; :\; \n ([0,t] \times \{ \zeta (\o )=\infty \}  ) 
\geq 1 \right\}  \\
\eta^* =\inf \left\{ t\geq 0 \; :\; \n^* ([0,t]\times \{ \zeta (\o )=\infty \}  ) \geq 1 \right\}  
\end{array} \right. \; ,$$
($\zeta (\omega ) $ being the lifetime 
of the path $\omega $). The random variables $\eta $ and $\eta^* $ have the same law as 
resp. $L_{\infty } $ and $L^*_{\infty }$, that are exponentially distributed or
infinite a.s. .

  In Section 3, Theorem
\ref{revexc} provides a decomposition of the law of the excursion under the supremum
reversed at its final jump: More precisely, we decompose the law of 
$(\o_{\zeta (\o)} - \o_{(\zeta (\o ) -s )-} ;\; 0\leq s \leq 
\zeta (\o) ) $ under $N^* (\cdot \cap \{ \o_{\zeta } >0 \}) $, in terms of the law of $X$ and its L\'evy measure. 
Theorems 
\ref{th} and \ref{revXlastsaut} 
give similar results for 
$(X_{\tau_x } -X_{(\tau_x -s )- } ;\; 0\leq s\leq \tau_x )$ 
under $\PP ( \; \cdot \; | X_{\tau_x } > x ) $
and $(X_{\sigma_x (t) } - X_{( \sigma_x (t) -s )-} ;\; 0\leq s \leq \sigma_x (t) )$
under $\PP ( \; \cdot \; | X_{\sigma_x (t) } > x ) $,  
where we have set for any $x, t>0 $ : 
$$ \tau_x =\inf 
\{ s\geq 0 :\; X_s > x \} \quad {\rm and }\quad \sigma_x (t) = \sup \{ s\in [0, t] :\; 
X_s \leq x \}\; . $$

  Williams in \cite{Wil74}, and many authors after him,
explored the connections between the Brownian motion, the three-dimensional Bessel process 
and the Brownian excursion (see for instance
Pitman \cite{Pit75} and Bismut \cite{Bis85}). Many of these identities 
in the Brownian case hold in a more general setting for totally 
asymetric L\'evy processes: 
see Bertoin \cite{Be92} for a generalized Pitman theorem for spectrally negatives L\'evy processes and 
Chaumont \cite{Ch}, \cite{Ch1} and \cite{Ch2} for Williams' theorems and Bismut's 
decomposition in the spectrally positive case. Let us mention that Chaumont has 
also explored the stable case in detail
in \cite{Ch97} and \cite{Ch}, providing several path-constructions and identities 
concerning the stable meander, the normalised excursion and the stable 
bridge. In these results the role of the three-dimensional Bessel process is played by the L\'evy process conditioned to stay positive. This process, 
denoted by $\Xpos $, has been introduced by Bertoin in a general setting (see \cite{Be93}). Bertoin's construction of $\Xpos $ 
is recalled in Section 2.2. We use it in combination with Theorem \ref{th} and \ref{revXlastsaut} to get in Section 4.1 
the generalized first 
Williams' decomposition theorem, then Bismut's decomposition of the excursion above the infimum 
in Section 4.2 and the second Williams' decomposition theorem in Section 4.3.

  Let us explain more precisely these results: For any $t>0 $, we define 
$U_t^* = X_{\locc_t }$ if $ L^*_{\infty }>t $ and $U^*_t = +\infty $ if not.
The process $(U^*_t ;\; t\geq 0 )$ is a subordinator 
(see Bertoin \cite{Be}, chapter 6) and its drift coefficient is denoted by
$d^* $. A classical result due to Kesten (see \cite{Kes69}) ensures  
that $\PP ( X_{\tau_x } = x )>0  $ iff $ d^* >0$. We assume that $d^* >0 $ and 
that $X$ does not drift to $-\infty $. Then, we can show that 
$\spos_x =\sup \{ s\geq 0 :\; \Xpos_s \leq x \} $ is finite a.s.. Theorem 
\ref{williamcont} show that 
$$ \PP ( \Xpos_{\spos_x } =x )= \PP (X_{\tau_x } =x ) $$ 
and that $ (x- X_{(\tau_x -s)-} ;\; 0\leq s\leq \tau_x  )$ 
under $\PP (\; \cdot \; \mid X_{\tau_x } =x ) $ has the same law as 
$(\Xpos_s ;\; 0\leq s \leq \spos_x ) $ under 
$ \PP (\; \cdot \; \mid \Xpos_{\spos_x } 
=x ) $.

  We also prove in Theorem \ref{truebismut} a path decomposition of 
the excursion above the infimum similar to Bismut's decomposition of the 
Brownian excursion: we show that
for any nonnegative measurable functionals $G$ and $D$ on the space of 
c\`adl\`ag paths with a finite lifetime
and for any nonnegative measurable function $f$,
\begin{multline*}
N \left( \int_0^{\zeta (\omega )} dt \; G\left( 
\o_s;\; 0\leq s\leq t \right) f(\o_t )
D\left( 
\o_{t+s}; \; 0\leq s\leq \zeta (\omega )-t \right) \right) \; = \\
\int_0^{+\infty } dx \; f(x)u^*(x) \EE \left[ G\left( 
\Xpos_s ;\; 0\leq s\leq \spos_x \right) 
\arrowvert 
\Xpos_{\spos_x } =x \; \right] 
\EE \left[ D\left( X_s ;\; 0\leq s\leq \tau_{-x }  \right) \right] \; ,
\end{multline*}
where $ \tau_{-x } = \inf \{ s\geq 0 :\; X_s < -x \} $ and where  
$u^*$ is the co-excessive version of the density of the potential 
measure associated with the subordinator $U^* $.

   Section 4.3 is devoted to the proof of Theorem \ref{williamdeux}
that can be seen as an analogue for general L\'evy 
processes of the second Williams' decomposition theorem that originally concerns the Brownian 
excursion split at its maximum. Let us describe our result: For any $x>0 $,
we set $\tpos_x = \inf \{ s\geq 0:\; \Xpos_s > x \} $. Proposition \ref{ind}
shows that 
$$\PP (\Xpos_{\tpos_x } =x )>0 \qquad {\rm iff} \qquad d^* >0 \; .$$ 
Let us denote by 
$\Xneg $ the process $X$ conditioned to stay negative (that is defined 
in Section 2.2); 
we write $\ggg (\o ) $ the instant when the 
excursion $\o $ attains its maximum. Theorem 
\ref{williamdeux} shows the law of $\o_{\ggg (\o ) } $ under $N$ admits a density with 
respect to Lebesgue measure that we specify. Under $N(\cdot \mid \o_{\ggg (\o ) } =x )$, the processes 
$ ( \o_s ;\; 0\leq s\leq \ggg (\o ))$ and 
$(\o_{s +\ggg (\o ) }; \; 0\leq s\leq \zeta (\o )-\ggg (\o ) ) $ are mutually
independent.
 Furthermore,

- the process $ ( \o_s ;\; 0\leq s\leq \ggg (\o ))$ is distributed as 
$( \Xpos_s ;\; 0\leq s\leq \tpos_x ) $ under
$\PP (\; \cdot\;  \mid \Xpos_{\tpos_x  } =x ) $;

- the law of $(\o_{s +\ggg (\o ) }; \; 0\leq s\leq \zeta (\o )-\ggg (\o ) ) $
 is absolutely continuous 
with respect to the law of $(\Xneg_s ;\; 0\leq s\leq \tneg_{-x } )$  
(with an evident notation for $\tneg_{-x} $)
and the corresponding density has the form $\varphi (\Xneg_{\tneg_{-x }} ) $,
where the function $\varphi $ is specified.

  Let us mention that we provide two other path decompositions that concern 
the excursion above the infimum (Theorem \ref{revexccont}) and the process
$(\Xpos_s ; \; 0\leq s\leq \spos_x )$ when $\Xpos_{\spos_x } >x $ (Theorem \ref{williamsaut}).

\section{Preliminary results.}

\subsection{Notation and basic assumptions.}

In this section we state our notation  
and the assumptions made at different 
stages of the paper. We also
recall fondamental results of
fluctuation theory that are our starting-point and we give some 
simple facts
concerning excursion theory applied to L\'evy processes, that is 
the main tool we use.

We begin with some notations concerning the canonical space. Let 
$\OO $ 
be the 
space of right-continuous functions with left limits from $(0, +\infty )$ to
$\RR $ (the so-called 
c\`adl\`ag functions space) endowed with the Skorokhod's topology. Let  
$\f $ stand for its Borel $\sigma $-algebra. For any path $\o $ in $\OO $
we define its lifetime $\zeta (\o ) $ by $\inf \{ t\geq 0 \; :\; \o (s)=
\o (t) \; , \; \forall s\geq t \} $, with the usual convention 
$\inf \emptyset =\infty $. For any time $t\geq 0 $, we denote the 
jump of $\o $ at $t$ by $\Delta \o (t)= \o (t)-\o(t- ) $; we also define 
the path respectively stopped at $t$, stopped just before $t$, reversed
at $t$ and reversed just before $t$, by 
\begin{eqnarray*}
&& \o (\cdot \wedge t)= \left(\o (s\wedge t) ;\; s\geq 0 \right), \qquad
\o(\cdot \wedge t- ) = ( \o(s\wedge t) -\Delta \o (t) 
\un_{[t, +\infty ) } (s) \; ;\; s \geq 0 ) , \\
&& \oo^t =\left( \o (t) -\o ( (t-s)- ); \; s\geq 0 \right), \qquad
\oo^{t- } =\oo^t -\Delta \o (t) ,
\end{eqnarray*}
\noi with the convention  $\o (s- )= \o (0) $ for any non-positive
real number $s$. When $\zeta (\o )$ is finite, $\oo^{\zeta (\o ) }$ is
well defined and simply denoted by $\oo $. 
We use {\bf a non-standard notation} for the shifted 
path at time $t$ defined by 
$$ \o \circ \theta_t = \left( \o (s+t ) -\o (t) \; ; \; s\geq 0 
\right) .$$ 

  For any $x>0 $, we denote by $\tau_x (\o ) $ and 
$\tau_{-x } (\o ) $ the first hitting time of respectively 
$(x , +\infty )$ and $(-\infty , -x )$:
$$ \tau_x (\o )= \inf \{s> 0 \; : \; \o (s) > x  \}, \qquad \tau_{-x }(\o )= \inf \{s> 0 \; : \; \o (s) <- x  \} \; ,$$
(with the usual convention $\inf \emptyset = + \infty $). 
For any time $t>0 $, we also denote by $\sigma_x (t,\o )$ and $\sigma_{-x} (t,\o ) $
the last passage time in  respectively $(-\infty , x]$ and 
$[-x , +\infty )$ on the time interval $[0,t]$ :
$$ \sigma_x (t,\o ) = \sup \{  0\leq s \leq t  \ \; : \; \o (s ) \leq x \} 
, \qquad \sigma_{-x} (t,\o ) = \sup \{  0\leq s \leq t  \ \; : \; \o (s ) 
\geq -x \} \; ,$$
(with the convention $\sup \emptyset = +\infty $). We write 
$\sigma_x (\o) = 
\lim_{t \rightarrow +\infty} \sigma_x (t, \o ) $, the limit being taken
in $[0, +\infty ] $. Next, we denote respectively by $\underline{g}_t (\o )$ and  
$\overline{g}_t (\o )$, the last infimum
time and the last supremum time of $\o $ before $t$:
$$\underline{g}_t (\o ) =\sup \{ s\in [0, t)  \; :\; \inf_{[0,t]} \o = 
\o (s- )\wedge \o (s) \}  $$
and 
$$ \overline{g}_t (\o ) =\sup \{ 
s\in [0, t)  \; :\;
\sup_{ [0,t] } \o =\o (s- )\vee \o (s ) \} . $$
\noi We also write $\underline{g} (\o ) =\lim_{ t\rightarrow +\infty }
\underline{g}_t (\o ) $ and 
$\overline{g} (\o ) =\lim_{ t\rightarrow +\infty }
\overline{g}_t (\o ) $ (note that these quantities may be infinite).

\vspace{7mm}

We denote by $X$ the canonical process on $\OO $: $X_t (\o )=\o (t) $ and 
we consider the probability measure $\PP $ on $(\OO , \f )$ under which $X$ 
is a L\'evy process started at $0$, with characteristic 
exponent $\psi $:

$$\EE \left[ e^{i\lambda X_t } \right] = e^{- t\psi (\lambda ) } \; , \quad 
t\geq 0 \; , \; \lambda \in \RR \; . $$

\noi By the L\'evy-Khintchine theorem, $\psi $ has the form
$$ \psi (\l )= ia\l + b\l^2 +\int  \pi (dr) \left( 1- e^{i\l r} +i\l r
\un_{ \{ |r| <1 \} } \right) \; , \quad  \l \in \RR , $$
\noi where $a$ is a real number, $b$ is  non-negative and the L\'evy measure $\pi $ 
is a Radon measure on $ \RR $ not charging $ 0 $, which  satisfies 
$$\int  \pi(dr) 
\left( 1\wedge |r|^2 \right)\; < \; +\infty . $$
\noi If  $J= \{ s\geq 0 \; :\; \Delta X_s \neq 0 \} $,  then
the point measure $\n (dsdr)=
\sum_{s\in J} \delta_{ (s, \Delta X_s )} $ is a Poisson measure with 
intensity $ds \pi (dr) $.

Let us recall some path-properties of L\'evy processes. For any $t\geq 0 $, we have 
$$ \widehat{X}^t \overset{(law)}{=} \left( X_s ;\; 0\leq s\leq t\right) \; .$$
(see Bertoin \cite{Be}). This identity is refered as the ``duality property''.

   In the whole paper (Section 3 excepted), we 
make the following assumption:

\vspace{5mm}

\noi {\bf Assumption (A) }: the point $0 $ is regular for $(0, +\infty )$ and  for $(-\infty ,0) $ 
with respect to $X$. 

\vspace{5mm}

\noi (In particular $X$ cannot be a subordinator or a compound Poisson process.) As a consequence of {\bf (A)}, we recall the following result
(see Millar \cite{Mil77}): For any $t\geq 0 $, the L\'evy process
$X$ reaches its infimum (resp. supremum) on $[0, t]$ at a unique instant 
that must be $\gi_t (X)$ (resp. $\gs_t (X)$).

  For every $t\geq 0$, we write 
$$S_t =\sup_{s\in  [0,t ]} X_s \; , \qquad I_t =\inf_{s\in [0,t ]} X_s .$$
\noi It is well-known that $X-S$ and $X-I$ are strong Markov processes (see Bertoin \cite{Be}, chapter 6). Assumption {\bf (A)} implies that $0 $ is regular for itself with respect to both these processes. Rogers has shown in 
\cite{Ro84} that this implies  
\begin{equation}
\label{rog1}
\PP \left( \exists t \in (0, +\infty ) \; :\; X_{t- } = I_{t-} < X_t \right) = 0
\end{equation}
and 
\begin{equation}
\label{rog2}
\PP \left( \exists t \in (0, +\infty )\; :\; X_{t- } = S_{t-} < S_t \right) = 0
\end{equation}
Let us recall briefly the proof: we only need to show 
for any $\epsilon >0 $
$$ \EE \left[ \sum_{s\in J } \un_{\{ X_{s-} = I_{s-} \} } 
\un_{ (\epsilon , +\infty ) } (\Delta X_s ) \right] \; =\; 0 \; .$$
Apply the compensation formula (see Bertoin \cite{Be} p. 7) to get 
$$ \EE \left[ \sum_{s\in J } \un_{\{ X_{s-} = I_{s-} \} } 
\un_{ (\epsilon , +\infty ) } (\Delta X_s ) \right] =
\pi ( (\epsilon , +\infty )) \int_0^{+\infty } ds \; \PP (X_s= I_s ) \; .$$
But the duality property implies for any $s>0 $, $\PP ( X_s = I_s ) =
\PP ( S_s = 0 ) = 0 $, because $0$ is regular for $(0, +\infty )$.
A similar argument proves (\ref{rog2}). \cqfd 

\vspace{7mm}

We denote the local times of $X-I$ and 
$X-S$ at the level $0$ by $(L_t )_{t\geq 0}$ and $( L^*_t )_{t\geq 0}$. They are uniquely determined  up to a multiplicative 
constant specified in a forthcoming lemma. The limit in $[0, +\infty ]$ of $L_t $ (resp. $L^*_t $) when $t$ goes to infinity is denoted by $L_{\infty }$
(resp. $L^*_{\infty } $). The quantity $L_{\infty } $ (resp. $L^*_{\infty } $)
is a.s. finite or a.s. infinite according as $X$ drifts or not 
to $+\infty $ 
(resp. $-\infty $). If $L_{\infty } $ (resp. $L^*_{\infty }$) is finite a.s., then it is exponentially distributed with parameter
denoted by $p$ (resp. $p^* $).

Equation (\ref{rog1}) and the dual result 
$ \PP ( \exists t \in (0, +\infty ) \; : \; X_{t- } = S_{t-} > X_t ) =0 $ imply 
that $\PP $-a.s. the sets 
$ \{ s\geq 0\; :\; X_s > I_s \} $ and $\{ s\geq 0 \; :\; X_s <S_s \} $  
are open sets (we have denoted $(g_i , d_i )\; , \; i\in \i$ and $(g_j , d_j )\; , \; j\in \i^* $
their respective connected components). Let $\mesu $ denote the Lebesgue 
measure on $\RR $. The duality property and assumption {\bf (A)} imply that 
$\mesu ( s\geq 0 \; :\; X_s =S_s )=\mesu (s\geq 0 \; :\; X_s =I_s )=0 $. Thus,
\begin{equation}
\label{nodrift}
\PP-{\rm a.s.} \qquad \mesu \left( \RR \backslash \bigcup_{i\in \i} (g_i , d_i ) \right)\; =\; 
\mesu \left( \RR \backslash \bigcup_{j\in \i^*} (g_j , d_j ) \right)\; 
=\; 0\; .
\end{equation}



Let $N$ and $N^*$ be the excursion measures of $X$ above its infimum and under its supremum as defined in the first section. 
Observe that as soon as the 
L\'evy measure $\pi $ charges $(-\infty, 0)$ (resp. $(0, +\infty )$), 
the set of excursions $\o $ ending with a negative jump (resp. positive jump)
has a positive $N$-measure (resp. $N^* $-measure). But thanks to 
(\ref{rog1}) and the dual result, we see that excursions above the
infimum and under the supremum leave $0$ continuously.

\vspace{7mm}

 Let $(\loc_t )_{t\geq 0}$ and $(\locc_t )_{ t\geq 0}$ be the right-continuous
inverses of $L$ and $L^* $: 
$$ L^{-1 }_t =\inf \{ s\geq 0 \; :\; L_s >t \} ,\qquad 
 L^{*-1 }_t =\inf \{ s\geq 0 \; :\; L^*_s >t \} \; ,$$
(with the convention $\inf \emptyset = \infty $). 
Recall that $\PP$-a.s.  
\begin{equation}
\label{intervalles}
\bigcup_{s\geq 0 } ( L^{-1}_{s-} , L^{-1}_{s} ) = \bigcup_{i\in \i} (g_i , d_i )  \qquad {\rm and } \qquad 
\bigcup_{s\geq 0 } ( L^{*-1}_{s-} , L^{*-1}_{s } ) = 
\bigcup_{j\in \i^*} (g_j , d_j ) \; .
\end{equation}

  For any $t\geq 0 $ we define  $U_t =-X_{\loc_t }$ if $ L_{\infty }>t $
and $U_t = +\infty $ if not. In a similar way, we define 
$U_t^* = X_{\locc_t }$ if $ L^*_{\infty }>t $
and $U^*_t = +\infty $ if not. The processes 
$(\loc , U )$ and $(\locc , U^* )$ are called the ladder processes. They are two-dimensional
subordinators killed at respective rates $p $ and $p^* $; their bivariate Laplace exponents are denoted by 
$$\kappa (\alpha , \beta )=
-\log \EE \left[ \exp (-\alpha \loc_1 -\beta U_1 ) \right] \quad {\rm and } \quad 
\kappa^* (\alpha , \beta )= - \log \EE \left[ \exp (-\alpha \locc_1 -\beta U^*_1 ) \right] $$  
\noi (see Bertoin \cite{Be}, chapter 6 for a detailed account). Next, we define 
the two potential measures $\U $ and $\U^* $ associated with $U $ and $U^* $:
$$\left\{ 
\begin{array}{ll}
\int_{ \RR  } \U (dx )f(x) = \EE \left[ \int_0^{L_{\infty }} dv \;  f(U_v ) \right] \\
\int_{ \RR } \U^* (dx )f(x) = \EE \left[ \int_0^{L^*_{\infty }} du \;  f(U^*_u ) \right] \; .
\end{array} \right . $$ 
\noi
Let $d^*$ be the drift coefficient of the subordinator $U^* $ :$\; d^* = \lim_{\beta \rightarrow +\infty } \kappa^* (0,\beta ) / \beta  $.
We recall the following result, due to Kesten \cite{Kes69} (see also Bertoin \cite{Be}, chapter 3, Theorem 5): assume that $d^* $ is positive,
and let $u^* \; : \; (-\infty , +\infty )\; \la [0, +\infty  ) \; $ be the co-excessive version of the density of $\U^* $. Then $u^* $  is 
continuous and positive on $(0,+\infty )$ , $\; u^* ( 0^+ )=1/d^* $, and 
$$\PP (X_{\tau_x } =x )= d^* u^* (x) \quad , \quad  x>0 \; , $$
\noi  where for convenience we write $\tau_x $ instead of $\tau_x (X) $. We prove the following simple lemma that will be used in Section 4.
\begin{lemma}
\label{pseudobismut}
Assume {\bf (A)}  and suppose that $d^* $ is positive. Then for any nonnegative 
measurable functional $F$ on $\OO $,
$$ \EE \left[  \int_0^{L^*_{\infty } } du F( X_{\cdot \wedge \locc_u } ) \right] \; =\; \int_0^{+\infty } dx \; u^* (x ) \EE \left[  F( X_{\cdot \wedge \tau_x } ) 
\; |\; X_{\tau_x }=x \; \right] \; . $$
\end{lemma}
\dem 
Set $A= \int_0^{+\infty } dx \; \EE \left[ F\left( X_{\cdot \wedge \tau_x } \right) \; ;\; X_{\tau_x } =x \right] $ and for any positive number 
$x$ define $H_x = L^*_{\tau_x } \; $ . Thanks to {\bf (A)}, we check path by path that
$$\PP-{\rm a.s.} \qquad {\rm on}\qquad \{ \tau_x < \infty \} \; ,\qquad 
\locc_{H_x } =\tau_x  \; .$$
\noi Thus
$$ A\; =\; \int_0^{+\infty } dx \; \EE \left[ F\left( X_{\cdot \wedge \locc_{H_x } }  \right) \; ;\; X_{\tau_x }=x \right] \; . $$
\noi Denote by $C$ the random set $\{ x\geq 0 \; : \; X_{\tau_x } =x \; \} $ and define the measures $\mu $ and $ \nu $ by
$$\left\{
\begin{array}{ll}
\mu (dx ) = \un_{C} (x) \; \mesu (dx) \; , \\
\int \nu (du) f(u) = \int \mu (dx ) f(H_x) \; .
\end{array} \right.  $$
Then,
\be
\label{pseudobismut1}
A \; =\; \EE \left[ \int_{[0, S_{\infty } )}  
\mu (dx ) F\left( X_{\cdot \wedge \locc_{H_x } } \right) \right] \; = \; \EE \left[ \int_{[0, L^*_{\infty }) } \nu (du ) 
F\left( X_{\cdot \wedge \locc_u } \right) \right] \; .
\end{equation}
 For any positive number $a$, 
\begin{eqnarray*}
\nu \left( [0,a]\right) & = & \mesu \left( \{ x\geq 0 \; : \; X_{\tau_x } =x \; ;\; H_x \leq a \; \} \right) \\
&= & \mesu \left( \{ x\geq 0 \; :\; X_{\tau_x } =x \; ; \; \tau_x \leq \locc_a \; \} \right) = \mesu \left( C\cap [0,U^*_a ] \right) \; .
\end{eqnarray*}
\noi Let us first consider the case $U^*_a <\infty $:
If there exists some $s$ in $[0,a ]$ such that $x\in (U^*_{s-} , U^*_s ) $, then, 
$\tau_x = L^{*-1 }_{s} $ and $ U^*_s=X_{\tau_x }>x $. Thus,
$$ \bigcup_{0\leq s\leq a}(U^*_{s-} , U^*_s ) \subset C^c \cap [0,U^*_a ] \; .$$
\noi Let $x$ be in $ C^c \cap [0,U^*_a ] $. Then, $X_{\tau_x } > x $. 
By (\ref{rog2}), it follows that $\tau_x $ must be the end-point of some excursion interval of 
$S-X$ above $0$ that is included in $[0, L^{*-1 }_a ] $. 
Then, by (\ref{intervalles}) there exists some $s$ in $[0,a]$ such that 
$$L^{*-1}_{s-} < L^{*-1 }_s =\tau_x \qquad {\rm and }\qquad S_{\tau_x -}=U^*_{s-} \leq x 
< U^*_s =X_{\tau_x } \; .$$
Hence,
$$ C^c \cap [0, U^*_a ] \subset  \bigcup_{0\leq s\leq a}[U^*_{s-} , U^*_s ) \; .$$
\noi By combining this with the previous inclusion we get 
$\mesu ( C^c \cap [0,U^*_a ] )= \sum_{0\leq s\leq a} \Delta U^*_s $.
But the Lévy-It\^o representation of $U^*$ guarantees that $\PP$-a.s.
$$ U^*_a = d^* a +\sum_{0\leq s\leq a} \Delta U^*_s \; ,\qquad 0\leq a < 
L^*_{\infty} \; .$$
\noi Then, $\PP$-a.s. for every $a$ in $[0, L^*_{\infty })\; $, 
$ \; \nu ([0, a] )=d^* a \; $. Next, observe that 
for any $a>L^*_{\infty }\; $,
$\; \nu ([0,a])= \mesu (C) = d^* L^*_{\infty } $. Thus, $\PP$-a.s.
$$ \nu (dx ) \; =\; d^* \un_{ [0, L^*_{\infty } ) } (x) m(dx ) $$
\noi The desired result follows from (\ref{pseudobismut1}) and the identity $\PP (X_{\tau_x }=x )=d^* u^* (x) $. \cqfd

\vspace{7mm}

  Let us introduce some notations:
for any positive time $t $ and any path $\o $, we denote the pre-infimum 
and the post-infimum path on the interval $[0,t]$ by:
$$\left\{ 
\begin{array}{ll}
\underleftarrow{\o^{t} } = \o ( \cdot \wedge \gi_t (\o ) ) \; , \\
\underrightarrow{\o^{t} } = (\o\circ \theta_{\gi_t (\o ) }) 
(\cdot \wedge (t-\gi_t (\o ) ) \; )  \; .
\end{array} \right. $$
\noi  We also denote the pre-supremum and post-supremum processes on $[0, t]$ by $\overleftarrow{\o }^t $ and $ \overrightarrow{\o }^t$:
$$\left\{ 
\begin{array}{ll}
\overleftarrow{\o }^t = \o ( \cdot \wedge \gs_t (\o ) ) \; , \\
\overrightarrow{\o }^t = (\o\circ \theta_{\gs_t (\o ) }) 
(\cdot \wedge (t-\gs_t (\o ) ) \; )  \; .
\end{array} \right. $$
We often use the following lemma in Section 4:
\begin{lemma}
\label{independance}
Assume {\bf (A)}. Let $T$ be independent of $X$ and exponentially 
distributed with parameter $\alpha >0 $. Then, $ \XgT $and $\XdT $ are mutually independent and the following identities hold for any nonnegative measurable 
functional $F$ on $\OO $:

\begin{description}
\item
$ (i) \qquad \EE \left[ F(\XgT ) \right] = \kappa (\alpha , 0) \EE \left[ 
\int_0^{L_{\infty } } dv e^{-\alpha L^{-1}_v } F( X_{\cdot \wedge L^{-1}_v } ) \right]
\; . $
 \item
$(ii) \qquad \EE \left[ F(\XdT ) \right] = \frac{1}{\kappa (\alpha , 0) }
N\left( \int_0^{\zeta } ds \alpha e^{-\alpha s }
F( \o_{ \cdot \wedge s} ) \right) \; . $
\end{description}
\end{lemma}
\dem
Let $G$ be any nonnegative measurable functional on $\OO $. We have 
$$ \EE \left[ F(\XgT ) G(\XdT ) \right]= \int_0^{\infty } dt \alpha 
e^{-\alpha t} \EE \left[ F(\Xg^t ) G(\Xd^t ) \right] \; . $$
\noi By (\ref{nodrift}) and by the definition of the excursions above the infimum, we have $\PP $-a.s.
$$ \int_0^{\infty } dt \alpha 
e^{-\alpha t} F(\Xg^t ) G(\Xd^t ) = \sum_{i\in \i } e^{-\alpha g_i }
F( X_{\cdot \wedge g_i } ) \int_0^{\zeta_i } ds \alpha e^{-\alpha s } 
G( \o^{i} (\cdot \wedge s) ) \; .$$
\noi Apply the compensation formula to get
$$\EE \left[ F(\XgT ) G(\XdT ) \right]= 
\EE \left[ 
\int_0^{L_{\infty } } dv e^{-\alpha L^{-1}_v } F( X_{\cdot \wedge L^{-1}_v } ) \right]
N\left( \int_0^{\zeta } ds \alpha e^{-\alpha s }
G( \o_{ \cdot \wedge s} ) \right) $$
and the following identities yield $(i) $ and $(ii) $:
$$  N\left( \int_0^{\zeta } ds \alpha e^{-\alpha s }
\right) = N\left( 1- e^{-\alpha \zeta }\right) \; =\; \kappa (\alpha , 0) $$
\noi and
$$ \EE \left[ 
\int_0^{L_{\infty } } dv e^{-\alpha L^{-1}_v } \right] = 
\frac{1}{\kappa (\alpha , 0) } \; .$$
\cqfd 

\vspace{7mm} 

We now specify the normalization of $L$ and $L^*$ thanks to the 
following proposition.

\begin{proposition}
\label{bismut}
Assume {\bf (A)}. Fix the normalization of $L $. Then, 
the normalization of $L^* $ can be chosen in order to have 
for any nonnegative measurable functional $F$ on $\OO $
$$ N \left( \int_0^{\zeta } ds \; F ( \wh \o^s ) \right) \; =\; \EE \left[ \int_0^{L^*_{\infty }} du \; F( X_{\cdot \wedge \locc_u } ) \right] \; $$ 
\noi and the dual identity
$$ N^* \left( \int_0^{\zeta } ds \; F ( \wh \o^s ) \right) \; =\; 
\EE \left[ \int_0^{L_{\infty }} dv \; F( X_{\cdot \wedge \loc_v } ) \right] \; . $$ 
\end{proposition}

\dem 
We denote by $\widehat{\XgT } $ the path $\XgT $ reversed at its lifetime 
$\underline{g}_T $. Observe that  
$$  \widehat{\XgT } = \widehat{X}^T \circ \theta_{ \overline{g}_T (
\widehat{X}^T ) } \; $$
(we use the fact that the minimum of $X$ over $[0, T]$ is attained $\PP$-a.s. at a unique time).
We also denote by $\widehat{\XdT }$ the path $\XdT $ reversed at its lifetime 
$ T-\underline{g}_T $. Similarly, we see that  
$$ \widehat{\XdT } = \widehat{X}^T_{\cdot \wedge 
\overline{g}_T ( \widehat{X}^T ) } \; .$$
The duality property implies that 
$$ \left( \widehat{\XgT } \; ,\; \widehat{\XdT } \right) \; 
\overset{(law)}{=} \; \left( \Xupd^T \; ,\; \Xupg^T \right) \; .$$
Let $G$ be any nonnegative measurable functional on $\OO $. Use Lemma
\ref{independance} to get:
$$ \EE \left[ F( \widehat{\XgT } ) G(\widehat{\XdT }) \right]= 
\EE \left[ 
\int_0^{L_{\infty } } dv e^{-\alpha L^{-1}_v} F( \widehat{X}^{L^{-1}_v } ) \right]
N\left( \int_0^{\zeta } ds \alpha e^{-\alpha s }
G( \oo^{s} ) \right) \; .$$
On the other hand, by replacing $X$ with $-X$, we see that Lemma \ref{independance} 
also implies 
$$ \EE \left[ F(\Xupd^T ) G(\Xupg^T ) \right]= 
\EE \left[ 
\int_0^{L^*_{\infty } } du e^{-\alpha L^{*-1}_u } G( X_{\cdot \wedge L^{*-1}_u } ) \right]
N^* \left( \int_0^{\zeta } ds \alpha e^{-\alpha s }
F( \o_{ \cdot \wedge s} ) \right) \; .$$
Thus, for any $\alpha >0 $
\begin{multline}
\label{bismut1}
\EE \left[ 
\int_0^{L^*_{\infty } } du e^{-\alpha L^{*-1}_u } G( X_{\cdot \wedge L^{*-1}_u } ) \right]
N^* \left( \int_0^{\zeta } ds e^{-\alpha s }
F( \o_{ \cdot \wedge s} ) \right)  \\
=\EE \left[ 
\int_0^{L_{\infty } } dv e^{-\alpha L^{-1}_v } F( \widehat{X}^{L^{-1}_v } ) \right]
N\left( \int_0^{\zeta } ds e^{-\alpha s }
G( \oo^{s} ) \right) \; .
\end{multline}
By letting $\alpha $ go to $0$, we see that the ratio 
$$ \displaystyle \frac{ N^* \left( \int_0^{\zeta } ds F( \o_{ \cdot \wedge s} ) \right) }{\EE \left[ 
\int_0^{L_{\infty } } dv F( \widehat{X}^{L^{-1}_v } ) \right]} $$
does not depend on $F$, provided it is well-defined (that is the denominator is positive and finite). Furthermore
this ratio coincides with 
$$ \frac{N \left( \int_0^{\zeta } ds G( \widehat{\o}^s ) \right) }{\EE \left[ 
\int_0^{L^*_{\infty } } du G( X_{\cdot \wedge L^{* -1}_u } ) \right]} $$
for any $G$ such that the denominator is positive and finite. We can choose the normalization of $N^* $, 
or equivalently $L^* $, so that both ratios are equal to $1$.  \cqfd 

\vspace{7mm}

In the spectrally positive case, the first identity of Proposition \ref{bismut} has been 
proved by Le Gall and Le Jan 
in \cite{LGLJ1} by a different method.

  Immediate applications of Proposition \ref{bismut} are the following identities due to Silverstein (see \cite{Sil80}), also mentioned in 
Rogers' paper \cite{Ro84} :
$$\left\{
\begin{array}{ll}
{\displaystyle N\left( \int_0^{\zeta } ds e^{-\alpha s-\beta \o_s } \right) }
\; =\; 1/\kappa^* (\alpha, \beta ) \\
{\displaystyle N^* \left( \int_0^{\zeta } ds e^{-\alpha s+\beta \o_s } \right) }
\; =\; 1/\kappa (\alpha, \beta ) \; . 
\end{array} \right . $$
\noi We can also derive from Proposition \ref{bismut} the Wiener-Hopf factorisation of the ladder exponents:
$$ \kappa (\alpha ,i\beta ) \kappa^* ( \alpha , -i\beta ) \; =\; \alpha + \psi (\beta ) .$$
Indeed, (\ref{nodrift}) gives the following decomposition:
$$ \int_0^{+\infty } dt \; e^{-\alpha t +i\beta X_t } \; =\; \sum_{i\in \i } \int_{g_i }^{d_i } dt  \; e^{-\alpha t +i\beta X_t } \; =\; 
\sum_{i\in \i } e^{-\alpha g_i +i\beta X_{g_i } }\int_0^{\zeta_i } ds \; e^{-\alpha s +i\beta \o_i (s) } \; \; . $$
\noi Taking the expectations and using the compensation formula, we get 
$$ 1/(\alpha + \psi (\beta ) ) \; =\; \EE \left[ \int_0^{L_{\infty }} du \; e^{-\alpha \loc_u +i\beta X_{\loc_u } } \right] 
N \left( \int_0^{\zeta }ds  e^{-\alpha s +i\beta \o_s } \right). $$
which yields the Wiener-Hopf factorization thanks to Proposition 
\ref{bismut}. 

\subsection{The Lévy process conditioned to stay positive or negative.}

We introduce now the process conditioned to stay positive, resp.  
negative, denoted by $\Xpos $, resp. $\Xneg $. 
Bertoin in \cite{Be93} provides 
a pathwise construction of $\Xpos $ and $\Xneg $ from concatenation
of the excursions of $X$ in $(0, +\infty )$, resp. $(-\infty , 0)$. Let us 
recall briefly this construction whose details can be found in 
\cite{Be93}, Section 3.

  Although Bertoin's construction holds in a general setting, 
{\bf we assume (A)}.
We denote by $(\f_t )_{t\geq 0} $ the natural filtration of $X$ completed 
with the $\PP $-null sets of $\f $. Then, $X$ is a semimartingale. Its 
continuous local martingale part is proportional to a standard Brownian 
motion and is independent of the non-continuous part. Let us denote 
by $\ell$ its semimartingale local time at $0$. We consider
$$ A^+_t = \int_0^t ds \un_{\{ X_s >0 \} } \qquad {\rm and }\qquad 
A^-_t = \int_0^t ds \un_{\{ X_s \leq 0 \} } \; .$$
Let us denote by $\alpha^+ $, resp. $\alpha^- $, the right-continuous 
inverse of $A^+ $, resp. $A^- $:
$$ \alpha^+_t =\inf \{ s\geq 0 \; :\; A^+_s > t \} \qquad
{\rm and} \qquad  \alpha^-_t =\inf \{ s\geq 0 \; :\; A^-_s > t \} \; ,$$
\noi (with the usual convention $\inf \emptyset =\infty $). Let $x$ be 
a real number. We denote its positive part, resp. negative part, by
$x_+ $, resp. $x_{-}$. We define a new  
process $\Xpos $ by 
$$ \Xpos_t = X_{\alpha^+_t } +\frac{1}{2} \ell_{\alpha^+_t} + \sum_{ 0<s\leq \alpha^+_t } 
\un_{ \{ X_s \leq 0 \} } (X_{s-} )_+ + \un_{\{ X_s >0 \} } (X_{s- })_{-} \quad {\rm if } \quad t< A^+_{\infty } $$
and by $\Xpos_t = +\infty $ if not. When $X$ has no 
Brownian part, $\Xpos $ can be viewed as the concatenation of 
the excursions of $X$ in $(0, +\infty )$. Similarly, we define $\Xneg $ by 
$$ \Xneg_t = 
X_{\alpha^-_t } -\frac{1}{2} \ell_{\alpha^-_t} -\sum_{ 0 <s\leq \alpha^-_t } 
\un_{ \{ X_s \leq 0 \} } (X_{s-} )_+ + \un_{\{ X_s >0 \} } (X_{s- })_{-} \quad {\rm if} \quad  t< A^-_{\infty } $$ 
and by $\Xneg_t = -\infty $ if not. 
The laws of $\Xpos $ and $\Xneg $ can be recovered by a harmonic transform:
Denote by $ q^+_t (x, dy ) $ and $ q^-_t (x, dy ) $ the semigroup of the 
Lévy process killed respectively in $(-\infty , 0] $ and $(0, +\infty )$.
One can show (see Silverstein in \cite{Sil80})
that the functions $\u^* ([0, x ]) $ and $\u ([0, x] ) $ are superharmonic
respectively for $q^+ $ and for $q^- $ and that the following 
kernels
$$ p^+_t (x, dy) = \frac{\u ([0,y ])}{\u ([0, x]} q^+_t (x, dy )
, \quad x> 0 $$ 
and
$$ p^-_t (x, dy ) = \ \frac{\u^* ([0, -y])}{\u^* ([0, -x])} q^-_t (x, dy )
, \quad x < 0 $$
define two sub-markovian semigroups. Bertoin has shown in \cite{Be93},
Theorem 3.4 that $\Xpos $ and $\Xneg $ are Markov processes 
started at $0$ with respective semigroups $p^+ $ and $p^- $. 
If $X$ does not drift to $-\infty $, resp. $+\infty $, then, 
$p^+ $, resp. $p^- $,
is markovian and $\Xpos $, resp. $\Xneg $, has an infinite lifetime. More precisely, if $X$ does not drift to $-\infty $, 
then, we can show that
\begin{equation}
\label{limite}
\Xpos_t <\infty \; , \quad t\geq 0 \qquad {\rm and } 
\qquad \lim_{t\rightarrow \infty } \Xpos_t = \infty \; .
\end{equation}
\noi {\bf Proof}. If $X$ does not drift to $-\infty $, it is easy to check that 
$\lim_{t\rightarrow \infty } A^+_t = \infty $, $\PP$-a.s. and then, 
$\Xpos_t <\infty  $, $t > 0 $. If $X$ drifts to $+\infty $, we have $\alpha^+_{t+A^+_{\sigma_0 }} = \sigma_0 +t;\ ,\quad t\geq 0 $, where 
$\sigma_0 = \sup \{ s\geq 0 :\; X_s \leq 0 \} < \infty $. Thus,
$$ \Xpos_{t+ A^+_{\sigma_0 } } \geq X_{\sigma_0 +t } \xrightarrow[t\rightarrow \infty]{\quad} +\infty \; .$$
If $X$ oscillates, we must consider to cases: Suppose first that $\pi \neq 0 $, then, 
$$  \lim_{t\rightarrow \infty } \Xpos_t \geq 
\sum_{s>0 } 
\un_{ \{ X_s \leq 0 \} } (X_{s-} )_+ + \un_{\{ X_s >0 \} } (X_{s- })_{-} \; 
= \infty \; .$$
If $\pi$ is null, then by assumption {\bf (A)} there is a Brownian component and  
$ \lim_{t\rightarrow \infty} \ell_t =\infty $ that yields  
the desired result. \cqfd 

\vspace{7mm}

  In particular cases, we recover ``classical''
definitions of the process conditioned to stay positive:

\noi {\bf -} In the Brownian case, 
$\u^* ([0, x ])= \u ((-x , 0]) =x $ and $p^+ $ is the semigroup of the three-dimensional Bessel process
started at $0$.

\noi {\bf -} In the spectrally positive case and the stable 
case, Chaumont has shown in \cite{Ch} and \cite{Ch2} that if the 
Lévy process does not drift to $-\infty $ and if $0$ is regular for 
$(0, +\infty )$, then, for any 
bounded $\f_t$ 
measurable functional $F$ that is continuous for the Skorokhod topology on 
$\OO $:
$$ \EE \left[ F( \Xpos ) \right] \; =\; 
\lim_{x \rightarrow 0} \; \lim_{T\rightarrow +\infty } \;
\EE_{x} \left[ F( X) \mid I_T \geq 0 \right] \; .$$

\noi {\bf -} In the spectrally negative case, Bertoin 
(see \cite{Be93}) gives 
another construction of $\Xpos $ that generalizes Pitman's 
theorem for Brownian motion (see Pitman \cite{Pit75}).

Let us denote 
by $\widehat{\Xg^t}$ the path $\Xg^t $ reversed at its lifetime
$\underline{g}_t $. We denote by $ (\Xneg_s ; \; 0\leq s < A^-_t ) $ (resp. $(\Xpos_s ; \; 0\leq s < A^+_t )$)
the process $\Xneg $ (resp. $\Xpos $) stopped at the random time $ A^-_t$ (resp. $ A^+_t$).
We need 
the following theorem due to Bertoin that links the process conditioned to 
stay positive with excursion theory:
\begin{theorem}
\label{bertoinresult} (Bertoin, \cite{Be93}, Theorem 3.1)
For every $t>0 $, the following identity holds 
$$ \left( \widehat{ \Xg^t }  \; ,\;  
\Xd^t \right) \; \overset{(law)}{=} \;
\left( (\Xneg_s)_{0\leq s < A^-_t }  
\; , \;  (\Xpos_s)_{0\leq s < A^+_t } \right) \; . $$
\end{theorem}

\vspace{7mm}

\rem If $X$ drifts to $+\infty $, then, the previous identity holds with  
$t=\infty $ as to say that $\Xpos $ has the same law as the post-infimum process 
(see Millar in \cite{Mil73}). 

\vspace{7mm}

   In Section 4, we use another identity that is proved in \cite{Be91} (see also \cite{LJ81}). 
From now on until the end of this section, we assume that $X$ does not drift to $-\infty$.  For any $t>0$ we set 
$$J_t = \inf_{s \in [t, \infty )} \Xpos_s  $$ 
and 
$$ \overline{d}_t =\inf \{ s>t  \; :\; S_t=X_s \} \; . $$
\begin{lemma}
\label{bertoinlemme} (Bertoin, \cite{Be91}, Lemme 4) 
$$ \left( S_{(\ddd_t +\ggg_t -t)-} -X_{(\ddd_t +\ggg_t -t)-} \; , \; S_{\ddd_t } \right)_{t \geq 0} 
 \; \overset{(law)}{=} \;
\left( \Xpos_t -J_t\; , \; J_t \right)_{ t\geq 0} \; .$$
\end{lemma}

\noi {\bf Proof:} Although Bertoin in Lemma 4 of \cite{Be91} only considers the case of $X \rightarrow +\infty $ (taking $\Xpos $ as the 
post-infimum process), the proof can be adapted when $X$ oscillates thanks to Theorem \ref{bertoinresult} and the arguments are exactly 
the same. \cqfd

\vspace{7mm}

The process $\Xpos -J$ is a strong Markov process and $0$ is a regular value. We denote by $K$ its local time at $0$ normalized 
in order it is distributed as $L^*$. Let us denote by $(g_i , d_i ) $ , 
$ i\in \i^{\uparrow} $ the excursion intervals of $\Xpos -J$ above $0$: 
$$ \{ s \geq 0 \; 
:\; \Xpos_s  > J_s \} = \bigcup_{i\in \ipos } (g_i , d_i ) \; , $$
We define
$$ w^i (s) = \Delta J_{g_i } + (\Xpos -J)_{(s+g_i )\wedge d_i } \; , \; s\geq 0 \; , \; i\in \ipos \; .$$
Then, Lemma \ref{bertoinlemme} implies that 
\begin{equation}
\label{poissonposi}
\n^{\uparrow} (dkdw) = \sum_{i\in \ipos } \delta_{(K_{g_i } , w^i )} 
\end{equation}
is a Poisson point process with intensity $dk \widehat{N}^* (dw)$, where $\widehat{N}^* $ is the law of $\widehat{\omega}^{\zeta}$
under $N^* (d\omega )$. We use this result in Section 4.

\section{Reversion formulas.}

Let $x$ and $t$ be two positive real numbers. 
We first decompose the law of $\xlast $
 on the event $\{ X_{\sigma_x (t)} >x \} $ in terms
of the law of $X$, the L\'evy measure $\pi $ and the function 
$\L $ that is defined on $(0, +\infty )\times (0, +\infty ) $ by $ \L (s,a ) \; =\; 
\PP (I_s \geq -a ) $. From classical fluctuation identities we have
$$ \int_{(0, +\infty )\times (0, +\infty ) } dsda 
\; \exp (-\l s -\mu a ) \L (s,a )  \; =\; \fr{\kappa (\l ,0 )}{\l \mu \kappa (\l , \mu ) }. $$
We also write $ \L (a ) $ for
the limit $\; \lim_{s\rightarrow  +\infty } \L (s,a) \; $ that is positive if and only if $X$ drifts to $ +\infty $ (or equivalently 
$p >0 \; $). To simplify notations, we write $\PP_r $ for the law of the L\'evy process started at $r $. 
We prove the first reversion formula:

\begin{theorem}
\label{revXlastsaut}
Assume that $\pi $ charges $(0, +\infty )$. Then for any positive numbers $x$ and $t$, and for any bounded measurable
functional $F$ on $\OO $,
\begin{multline*}
\EE \left[  F(\xlast ) ;\; \sigma_x (t) <\infty ;\; 
\; \xsigma > x \; \right]  \\
=\int_{(0,+\infty ) } \pi (dr)  \int_0^t du  \EE_r \left[  F( X_{\cdot \wedge u} ) \L (t-u , X_u  -x ) ;\; 
x < X_u \leq x +r  \right] \;  . 
\end{multline*}
\end{theorem}

\vspace{7mm}

\conss  (i) If $X$ drifts to $+\infty $, then $\L (a)> 0$ for any positive real number $a$. Thanks to Theorem
\ref{revXlastsaut} we get
\begin{multline}
\label{tpsinfty}
\EE \left[  F(\xlast ) \frac{\L (\xsigma -x)}{\L (t-\sigma_x (t) , \xsigma -x )} ;\; \sigma_x (t) <\infty ;\; 
\; \xsigma  > x \; \right] 
\\
= \int_{(0,+\infty) } \pi (dr)  
\int_0^{t} du \; \EE_r \left[  F( X_{\cdot \wedge u} ) \L ( X_u -x ) \; ;\; x< X_u \leq x+r \; \right] \;  . 
\end{multline}
Observe that $\PP$-a.s. $\sigma_x (t)= \sigma_x $ for all $t$ sufficiently large. 
Thus, $\PP$-a.s.
$$ F(\xlast ) \frac{\L (\xsigma -x)}{\L (t-\sigma_x (t) , \xsigma -x )} 
\un_{ \{  \sigma_x (t) <\infty ;\;\xsigma  > x \} }
\xrightarrow[t \rightarrow \infty ]{\; }
F(\wh X^{\sigma_x } )  \un_{ \{ X_{\sigma_x } > x \} } \; .$$
Since $\L (\xsigma -x) / \L (t-\sigma_x (t) , \xsigma -x ) $ is smaller than $1$, dominated convergence applies and 
we deduce from (\ref{tpsinfty}) that
\begin{multline*}
\EE \left[  F(\wh X^{\sigma_x } ) \; ;\; X_{\sigma_x } > x \; \right]  \\
= \int_{(0,+\infty ) } \pi (dr)  
\int_0^{+\infty } du \; 
\EE_r \left[  F( X_{\cdot \wedge u} ) \L ( X_u -x ) 
\; ;\; x< X_u \leq x +r\; \right] \; . 
\end{multline*}

\vspace{7mm}

\noi (ii) By the duality property applied in  the right side of Theorem \ref{revXlastsaut}, we see that 
under $\PP ( \cdot \mid X_{\sigma_x (t) }>x ) $
$$ X_{\cdot \wedge \tpssig- } \; \overset{(law)}{=}\; \wh X^{\tpssig- } \; .$$

\vspace{7mm}

\noi {\bf Proof of Theorem \ref{revXlastsaut}.}
Let $\epsilon $ be a positive real number and let $(\sigma_n )_{n\geq 0 } $ be the increasing sequence
of the jump times $\{ s \geq 0 \; :\; \Delta X_s > \epsilon \} $. Recall that $\sum_{s\in J}\delta_{(s,\Delta X_s )} $
is a Poisson measure with intensity $dl\pi (dr) $. Let $f$ 
be a bounded function on $\RR $. Consider the event  $ A_{\epsilon } =
\{ \sigma_x (t) <\infty ;\; \inf_{s\in [\sigma_x (t), t ] } X_s  > \epsilon +x \} $ and 
set 
$$ a (\epsilon )= \EE \left[ F( \widehat{X}^{\sigma_x (t) -} ) f( \Delta 
X_{\sigma_x (t) } ) \; ;\; A_{\epsilon} \right] \; .$$
Observe that 
$$ a(\epsilon ) = \sum_{n\geq 0 }
\EE \left[ F( \widehat{X}^{\sigma_n - } ) f( \Delta X_{\sigma_n } ) ;\; 
\sigma_n <t ; \; X_{\sigma_n -} < x ;\; X_{\sigma_n } + \inf_{[0, t-\sigma_n ] } X\circ 
\theta_{\sigma_n } > \epsilon +x \right] .$$
Apply the Markov property at $ \sigma_n $ in order to get
$$  a(\epsilon ) = \sum_{n\geq 0 }
\EE \left[ F( \widehat{X}^{\sigma_n - } ) f( \Delta X_{\sigma_n } ) 
\L ( t- \sigma_n , X_{\sigma_n } -\epsilon -x ) \; ;\; 
\sigma_n <t ; \; X_{\sigma_n -} < x   \right]  .$$
Then,
$$ a(\epsilon ) = \EE \left[ \sum_{ s\in J \; : \; s\leq t } 
\un_{\{ \Delta X_s > \eps ;\;  X_{s- } <x \} } F( \widehat{X}^{s- } ) f(\Delta X_s )
\L (t-s , X_{s-} +\Delta X_s -\epsilon -x ) 
\right] .$$
Apply the compensation formula to get:
$$  a(\epsilon )= \int_{(\epsilon , +\infty )} \pi (dr) f(r) \int_0^t du \EE 
\left[ F( \widehat{X}^{u} ) \L ( t-u, X_u +r -\epsilon -x ) ; \; 
x+\eps - r < X_u <x 
\right] $$
and by duality 
\begin{equation}
\label{findus}
a(\epsilon )= \int_{(\epsilon , +\infty )} \pi (dr) f(r) \int_0^t du \EE 
\left[ F( X_{\cdot \wedge u} ) \L ( t-u, X_u +r -\epsilon -x ) ; \;
x+\eps - r < X_u <x 
\right] \; . 
\end{equation}
Next, observe that 
$\PP $-a.s. 
$$ \lim_{\epsilon \rightarrow 0} 
\un_{ A_{\epsilon } } =\un_{ \{ \sigma_x (t)<\infty ;\; 
X_{\sigma_x (t)} > x  \} } \; $$
and complete the proof by letting $\epsilon $ go to $0$ and using dominated convergence in the 
left side of (\ref{findus}) and monotone convergence in the right side. \cqfd 

\vspace{7mm}

We get a similar result for the 
reversed path at $\tau_x $ on the event $\{ X_{\tau_x } >x \} $:
\begin{theorem} 
\label{th}
Assume that $\pi $ charges $(0, +\infty ) $. Then, for any 
positive real number $x$ and for any 
bounded measurable functional $F$ on $\OO $, 
\begin{multline*}
\EE \left[ F(\wh X^{\tau_x } ) ;\;  \tau_x <\infty ;\;  \; X_{\tau_x }>x  \right] \\
=\int_{ (0,+\infty ) } \pi (dr ) \int_0^{+\infty }du \;  
\EE_r \left[ F( X_{\cdot \wedge u } ) \; ;\; x< X_u \leq x + I_u \; \right] \; .
\end{multline*}
\end{theorem}

\vspace{7mm}

\rem
In the subordinator case we get immediately the well-known formula:

$$ \EE \left[  f(X_{\tau_x- }, X_{\tau_x } )\; ;\;  X_{\tau_x } >x \right] \; =\; 
\int_{[0, x] } V(da )
\int_{(x-a, +\infty )} \pi (dr ) f( a , a+ r ) \; ,$$
\noi where $V$ denote the potential measure associated with $X$.

\vspace{7mm}

\noi {\bf Proof of Theorem \ref{th}.}
Let $f$ be a bounded measurable function on $\RR $. Observe that
\begin{multline*}
\EE \left [ F( \widehat{X}^{\tau_x - } ) f( \Delta X_{\tau_x } ) ; \;
\tau_x <\infty ;\; X_{\tau_x } > x \right] 
\\ = \EE \left[ \sum_{s\geq 0 } 
\un_{ \{ S_{s- } \leq x \; ;\; \Delta X_s + X_{s- } > x \} } F(\widehat{X}^{s- } ) f(\Delta X_s ) 
\right] \; .
\end{multline*} 
Apply the compensation formula to get 
\begin{multline*}
\EE \left [ F( \widehat{X}^{\tau_x - } ) f( \Delta X_{\tau_x } ) ; \;
\tau_x <\infty ;\; X_{\tau_x } > x \right]  \\
=\int_{(0, +\infty )} \pi (dr ) f(r) \int_0^{+\infty } du \; 
\EE \left[ F(\widehat{X}^u ) ;\; S_u \leq x ;\; r+X_u >x \right] 
\end{multline*}
\noi and the result follows by the duality property. \cqfd

\vspace{7mm}

  Recall that if the L\'evy measure charges $(0, +\infty )$, 
the excursion under the supremum may end with a jump. We now give  
a decomposition of the law of the excursion under the supremum reversed at its final jump time.

\begin{theorem}
\label{revexc}
Assume {\bf (A)} and that $\pi $ charges $(0, +\infty )$. Then, for 
any nonnegative measurable functional $F$ on $\OO $,  
$$N^* \left( F(\wh \o^{\zeta } )\; ;\; \o_{\zeta } >0 \right) \; =\; \int_{(0, +\infty)} \pi (dr )\EE_r \left[ \int_0^{L_{\infty } }dv \; 
 F( X_{\cdot \wedge \loc_v } ) \un_{\{ X_{\loc_v } >0 \} } \right] .$$
\end{theorem}

\vspace{7mm}

\rem  In the spectrally positive case, if we take $L=-I$, then, 
$L^{-1}_x = \tau_{-x} $. Theorem \ref{revexc} shows that under
$N^* (\cdot \cap \{ \o_{\zeta } >0 \} )$,
 the law of $\o_{\zeta -}$ admits 
a density with respect to Lebesgue measure that is given by
$$ x \longrightarrow \un_{(-\infty, 0)} (x) \PP ( I_{\infty }< x ) \pi ( (-x , +\infty )) \; $$
\noi Furthermore, under $N^* (\cdot \mid \o_{\zeta -} =-x )$, the path $\widehat{\o }^{\zeta -}$
is distributed as $ X_{\cdot \wedge \tau_{-x} }$ under $\PP (\cdot \mid \tau_x <\infty ) $.  
This result has been used by Bertoin in \cite{Be92} and \cite{Be91}.

\vspace{7mm}

\noi {\bf Proof of Theorem \ref{revexc}.} For any nonnegative measurable function $f$ on $\RR $, we have the following decomposition:
$$ N^* \left( F(\oo^{\zeta -} ) f( \Delta \o_{\zeta } ) ;\; \o_{\zeta }>0 \right)
=N^* \left( \sum_{s\geq 0 } \un_{\{ \o_{s-} +\Delta \o_s >0 \} } F(\oo^{s- } ) f(\Delta \o_s ) \right) \; .$$
Observe that $\o $ under $N^* $ is markovian with the transition kernel of the L\'evy process killed in $[0, +\infty )$. 
We can apply the compensation formula to get 
$$  N^* \left( F(\oo^{\zeta -} ) f( \Delta \o_{\zeta } ) ;\; \o_{\zeta }>0 \right)
= N^* \left( \int_0^{\zeta } ds \; F( \oo^{s} ) \int_{ (0, +\infty )} \pi (dr ) f(r) \un_{ \{ \o_s +r >0 \} } \right) $$
that yields the theorem thanks to Proposition \ref{bismut}. \cqfd

\section{Applications.}

\subsection{First Williams' decomposition theorem.}

{\bf From now on until the end of the present article, we assume (A) 
and we suppose that $X$ does not drift to $-\infty $}
Let $x$ be a positive real number. Williams has 
shown in \cite{Wil74} that the standard real Brownian motion 
reversed at the first hitting time of $(x, +\infty )$ is 
distributed as the three-dimensinal Bessel process up to 
its last passage time at $x$. In this section, we extend Williams' result to 
general L\'evy processes, the role of the three-dimensional Bessel process
being played by the L\'evy process conditioned to stay positive.
In order to avoid cumbersome notation, we set
for any positive real numbers $x$ and $t$
\begin{eqnarray*}
\sd_x^{\scriptscriptstyle{t} } = \sigma_x (\Xd^{\scriptscriptstyle{t} } ) &=& 
\sup \{ s\in [0, t-\underline{g}_t ] : \; X_{s + \underline{g}_t } - I_t \leq x \} \; , \\
\spos_x = \sigma_x ( \Xpos ) &=& \sup \{ s\geq 0 :\; \Xpos_s \leq x \} \; , \\
\td_x^{\scriptscriptstyle{t}} = \tau_x (\Xd^{\scriptscriptstyle{t} } ) &=& \inf \{ s \in [0, t-\underline{g}_t ]:\;    X_{s + \underline{g}_t } -I_t > x \} 
\; ,\\
\tpos_x = \tau_x ( \Xpos ) &=& \inf \{ s\geq 0 :\; \Xpos_s >x \} \; ,
\end{eqnarray*}
(with $ \inf \emptyset = \sup \emptyset = +\infty $). Observe that $\sd^t_x $ may be infinite if 
$X_t -\inf_{[0, t]} X < x $. Under $\PP ( \cdot \cap \{  \sd^t_x < \infty \} ) $ we define 
$$ Y^{ t } \; =\; \Xd^t \circ \theta_{ \sd_x^{\scriptscriptstyle{t} } } $$
Similarly we denote 
$$ Y \; =\; \Xpos \circ \theta_{\spos_x } \; $$
that is well-defined thanks to (\ref{limite}) and our assumptions. Notice 
that $Y^t$ and $Y$ rely on $x$ although it does not appear in the notations. We recall that $\Yg$ and $\Yd $ are respectively 
the pre-infimum process and the post-infimum process of $Y$. The following theorems describe the law of the path $\Xpos $ reversed
at time $\spos_x $ : the first theorem concerns the case of a jump : $\Delta \Xpos_{\spos_x } >0 $; the second theorem deals 
with the process leaving continuously level $x$. 
 \begin{theorem} 
\label{williamsaut}
Assume that $\pi $ 
charges $(0, +\infty )$. Let $x>0 $.
\begin{description}
\item
(i) $ \left( \widehat{\Xpos }^{\spos_x }\; ,\; \Yg \right) $ under $ \PP ( \cdot \mid \Xpos_{\spos_x } > x ) $ 
$\overset{(law)}{=} \; \left( \; \Delta X_{\tau_x } + X_{\cdot \wedge 
\gs_{\tau_x} } \; ,\; \widehat{X}^{\tau_x -}_{\cdot \wedge (\tau_x - 
\gs_{\tau_x} )} \; \right) $ under $\PP ( \cdot \mid X_{\tau_x } > x ) $. 
\item
(ii) Under  $ \PP ( \cdot \mid \Xpos_{\spos_x } > x ) $, $\Yd $ is independent of 
$ (\Xpos_{\cdot \wedge \spos_x }\; ,\; \Yg )$ and distributed as $\Xpos $.
\end{description}
\end{theorem}

\begin{theorem} 
\label{williamcont}
Assume that $d^* >0 $. Then, for any $x>0 $,
\begin{description}
\item
(i) $ \PP ( \Xpos_{\spos_x } =x )= \PP (X_{\tau_x } =x )= d^* u^*(x) $ and 
$\widehat{X}^{\tau_x } $ under
$\PP (\cdot \mid X_{\tau_x } =x ) $ is distributed as 
$\Xpos_{\cdot \wedge \spos_x } $ under $\PP (\cdot \mid \Xpos_{\spos_x } 
=x )$.

\item
(ii) Under  $ \PP ( \cdot \mid \Xpos_{\spos_x } =x ) $, $\Xpos \circ \theta_{\spos_x } $ and 
$\Xpos_{\cdot \wedge \spos_x } $ are mutually independent and 
$\Xpos \circ \theta_{\spos_x } $ is distributed as $\Xpos $. 
\end{description}
\end{theorem}

\vspace{7mm}

\rem We assume that $\pi $ charges $(-\infty , 0 )$. Then, the 
excursion under the infimum may end with a negative jump. The dual form of the reversion formula of 
Theorem \ref{revexc} gives
\begin{equation}
\label{excucur}
 N\left( F(\o_{\cdot \wedge \zeta- } ) \; ; \; \o_{\zeta } <0 \right) 
 = \int_{(-\infty , 0 ) } \pi (dr ) \EE \left[ \int_0^{L^*_{\infty }} du \un_{\{ X_{L^{*-1}_u } <-r \} } 
F( \widehat{X}^{L^{*-1}_u } ) \right] \; .
\end{equation}
We assume moreover that $d^* >0 $. By 
Lemma \ref{pseudobismut}, it follows that under $N ( \cdot \cap \{ \o_{\zeta } <0 \}) $, $\o_{\zeta -} $ admits a density with respect to 
Lebesgue measure given by
$$ N\left( \o_{\zeta -} \in dx \; ;\; \o_{\zeta }<0 \right) = u^* (x) \pi ( (-\infty , -x)) {\rm m} (dx ) \; .$$ 
By combining (\ref{excucur}) with Lemma \ref{pseudobismut}, Theorem \ref{williamcont} implies 
that for any $x>0 $, 
$$ \Xpos_{\cdot \wedge \spos_x } \quad {\rm under} \quad \PP (\cdot \mid \Xpos_{\spos_x } =x )
\overset{(law)}{=} \o_{\cdot \wedge \zeta -} \quad {\rm  under} \quad
N( \cdot \mid \o_{\zeta -} =x ) \; .$$
This result is due to Chaumont in 
the stable case (see 
\cite{Ch}).

\vspace{7mm}

\noi {\bf Proof of Theorem \ref{williamsaut}.} First observe that $\tau_x (J)=\spos_x $ a.s., then 
$$ \spos_x \; \overset{(law)}{=}\; \inf \{ t> 0 \; : \; S_{\ddd_t } >x \} = \ggg_{\tau_x } \; .$$
Set 
$$ \gamma = \spos_x + \gg (Y) = \inf \{ t>\spos_x  \; : \; \Xpos_t = J_t  \} \; .$$
Then, 
\begin{equation}
\label{taugamma}
\gamma \; \overset{(law)}{=} \;  \tau_x \; . 
\end{equation}
We define the functional $\Sigma $ by 
$$ \Sigma_t (X) = \left( S_{(\ddd_t + \ggg_t -t)- } - X_{(\ddd_t + \ggg_t -t)- }\; , \;   S_{\ddd_t }  \right) \;\; , \; t\geq 0 \; .$$
Deduce from (\ref{taugamma}) and from the fact that $ \ddd_{\tau_x } = \tau_x $ that
\begin{equation}
\label{independancesigma}
\left( (\Xpos -J)_{\gamma +t} \; , \; J_{\gamma +t} -J_{\gamma} \right)_{t\geq 0}  
\overset{(law)}{=} \left( \Sigma_t (X\circ \theta_{\tau_x } ) \right)_{t\geq 0} \; . 
\end{equation}
Since 
$$ \Yd_t = (\Xpos -J)_{\gamma +t} + (J_{\gamma +t} -J_{\gamma}) $$
we deduce from the Markov property applied at $\tau_x $ in the right member of 
(\ref{independancesigma}) that $\Yd $ is independent of $\Xpos_{\cdot \wedge \gamma}$ and that $\Yd $ has the same distribution as $\Xpos $, 
which proves Theorem \ref{williamsaut} (ii) and also Theorem \ref{williamcont} (ii) because 
$$ \Xpos \circ \theta_{\spos_x } = \Yd \quad {\rm on } \quad \{ \Xpos_{\spos_x } =x \} \; . $$
\noi Next, we denote by $K^{-1 } $ the right-continuous inverse of $K$:
$$ K^{-1}_u = \inf \{ t\geq 0 \; :\; K_t >u \} \quad u\geq 0 .$$
We need the following lemma:
\begin{lemma}
\label{dualityLJ} For any x>0 
$$ \Xpos_{\cdot \wedge K^{-1}_x } \; \overset{(law)}{=} \; \widehat{X}^{L^{*-1}_x } \; .$$
\end{lemma}
\noi {\bf Proof:} We first index the excursions of $\Xpos -J$ above $0$ 
by the corresponding local time: for any $t\geq 0$, we set 
$$ e^{\uparrow}_t = \left( (\Xpos -J)_{(K^{-1}_{t-} +s )\wedge K^{-1}_t} \right)_{s\geq 0} $$
(note that $e^{\uparrow}_t =0$ if $\Delta K^{-1 }_t =0 $). We do the same thing for the excursions of $X$ under its supremum and we set 
$$ e_t = \left( (X-S)_{(L^{*-1}_{t-} +s )\wedge L^{*-1}_t} \right)_{s\geq 0} \; .$$
We also write $U^{\uparrow}_t = \Xpos_{K^{-1}_t }$. We deduce from Lemma \ref{bertoinlemme} that $(K^{-1} , U^{\uparrow})$ is a subordinator 
with the same distribution as $(L^{*-1} , U^*)$. There exists a measurable functional $F$ such that 
$$ F\left( (e^{\uparrow}_t , U^{\uparrow}_t )_{0\leq t\leq x} \right) = \Xpos_{\cdot \wedge K^{-1}_x } \; , \; a. s. \; .$$
Let us explain more precisely how to recover $\Xpos $ from the $e^{\uparrow}_t $ and $U^{\uparrow}_t $, $t\geq 0$: For any 
$s\in [0, K^{-1}_x ]$, we 
define 
$$ g(s) = \sup \{ u\in [0, s) \; :\; \Xpos_u =J_u \} \; .$$
The L\'evy-Ito decomposition for the subordinator $K^{-1} $ implies that 
$$ K_s = \sup \left\{ a\in [0, x] \; :\; 
d^* a + \sum_{ t<a } \zeta (e^{\uparrow}_t ) \leq s \right\} $$ 
and 
$$ g(s) = d^* K_s + \sum_{t< K_s } 
\zeta (e^{\uparrow }_t )  \; .$$
Then, 
$$ \Xpos_s = e^{\uparrow}_{K_s } (s-g(s) ) + U^{\uparrow}_{K_s } \; .$$
In order to simplify notations, we set for any $t \geq 0$:
$$ \widehat{e}_t = \widehat{e}^{\zeta (e_t)- }_t \; .$$
Lemma \ref{bertoinlemme} implies 
\begin{equation}
\label{premiereidentite}
\left( e^{\uparrow}_t , U^{\uparrow}_t \right)_{0\leq t\leq x}  \; \overset{(law)}{=} \;
\left( \widehat{e}_t , U^*_t \right)_{0\leq t\leq x} \; .
\end{equation}
Since $(e_t \; ; \; t \geq 0)$ is a Poisson process and $U^*$ a subordinator, a simple time-reversal argument show that
$$ \left( e_{x-t} , \widehat{U^*}^x_t \right)_{0\leq t\leq x}  \; \overset{(law)}{=} \; 
\left( e_t , U^*_t \right)_{0\leq t\leq x} \; . $$
Thus 
\begin{equation}
\label{deuxiemeidentite}
\left( e^{\uparrow}_t , U^{\uparrow}_t \right)_{0\leq t\leq x}  \; \overset{(law)}{=} \;
\left( \widehat{e}_{x-t} , \widehat{U^*}^x_t \right)_{0\leq t\leq x}
\end{equation}
Applying the L\'evy-Ito decomposition for the subordinator $L^{*-1}$ reversed at time $x$, it is easy to check that 
$$ F \left( ( \widehat{e}_{x-t} , \widehat{U^*}^x_t )_{0\leq t\leq x} \right) = \widehat{X}^{L^{*-1}_x } \; , \quad a.s. $$
and we conclude thanks to (\ref{deuxiemeidentite}). \cqfd

\vspace{4mm}

Let us prove now Theorem \ref{williamsaut} (i): Let $F$ and $G$ be two nonnegative measurable functionals and $f$ be a nonnegative 
measurable function. Set 
$$ \alpha = \E \left[ F\left( \Xpos_{\cdot \wedge \spos_x -} \right) f(\Delta \Xpos_{\spos_x } ) G( \Yg ) \; ; \; 
\Xpos_{\spos_x } > x \right] \; .$$
It is sufficient to show that 
\begin{equation}
\label{crucialpoint}
\alpha = \E \left[ F \left( \widehat{X}^{\ggg_{\tau_x } } \right) f( \Delta X_{\tau_x } ) G 
\left( \widehat{X}^{\tau_x -}_{\cdot \wedge (\tau_x - 
\gs_{\tau_x} )} \right) \; ; \; X_{\tau_x } >x \right] \; . 
\end{equation}
First observe that 
$$ \alpha = \E \left[ \sum_{ i\in \ipos } \un_{ \{ \Xpos_{g_i } \leq x < 
\Xpos_{g_i } + w^i (0) \} } F ( \Xpos_{\cdot \wedge g_i } ) f(w^i (0) )
G( w^i -w^i (0) ) \right] \; .$$
Then by \ref{poissonposi}, we get
$$ 
\alpha =\int_0^{\infty} du \E \left[ F( \Xpos_{\cdot \wedge K^{-1}_u } ) \widehat{N}^* \left( f(w (0)) G( w -w(0) );\; 
U^{\uparrow}_u \leq x < U^{\uparrow }_u + w(0)  \right) 
\right] \; .$$
The previous lemma implies that 
\begin{equation}
\label{excursioncruciale}
\alpha = \int_0^{\infty} du \E \left[ F( \widehat{X}^{L^{*-1}_u} ) N^* \left( f(\Delta \omega_{\zeta}  ) G( \widehat{\omega}^{\zeta -} )
;\; U^*_u \leq x < U^*_u + \omega (\zeta )   \right) \right] \; .
\end{equation}
But we have a.s.
$$ \un_{ \{ X_{\tau_x } >x \} } F \left( \widehat{X}^{ \ggg_{\tau_x } } \right) f ( \Delta X_{\tau_x } ) G 
\left( \widehat{X}^{ \tau_x -}_{\cdot \wedge ( \tau_x - \gs_{\tau_x} ) }   \right)  = \sum_{ i \in \i^*} 
\un_{ \{  X_{g_i } \leq x < X_{g_i } + \omega^i (\zeta_i) \} } 
F ( \widehat{X}^{g_i } ) f(\Delta \omega^i_{\zeta_i} ) 
G( \widehat{ \omega^i}^{ \zeta_i -} ) $$
Then, the compensation formula combined with (\ref{excursioncruciale}) achieve the proof of (\ref{crucialpoint}). \cqfd

\vspace{7mm}

\noi {\bf Proof of Theorem \ref{williamcont}:} We only need to show $(i)$. From Theorem \ref{williamsaut}, we deduce that
$$ \PP \left( \Xpos_{\spos_x } > x \right) = \PP \left( X_{\tau_x } >x \right) $$ 
So, we have for any $x>0 $:
\begin{equation}
\label{marginalbismut}
d^* u^* (x) = \P ( \Xpos_{\spos_x } =x  ) \; . 
\end{equation}
We need the following lemma:
\begin{lemma}
\label{balayageposi} Under the assumptions of Theorem \ref{williamcont}, we have for any nonnegative measurable functional $F$
$$ \E \left[ \int_0^{\infty} du \F ( \Xpos_{\cdot \wedge K^{-1}_u }) \right] = \int_0^{\infty} dx u^* (x) \E \left[ 
F( \Xpos_{\cdot \wedge \spos_x }) \; \mid \; \Xpos_{\spos_x } =x \right] \; .$$
\end{lemma}

\noi {\bf Proof:} 
We argue exactly as in Lemma \ref{pseudobismut} replacing, $X$ by $\Xpos $,  $L^*$ by $K$ and $\tau_x $ by $\spos_x $. \cqfd 

\vspace{7mm} 

Lemmas \ref{pseudobismut} and \ref{balayageposi} imply that for any nonnegative measurable functional $F$, the set
$$ \l_F = \left\{ x >0 \; :\; \E \left[ F( \Xpos_{\cdot \wedge \spos_x }); \Xpos_{\spos_x } =x \right] =
\E \left[ F( \widehat{X}^{\tau_x }) ; X_{\tau_x } =x \right] \; \right\} $$
is of full Lebesgue measure. We have to show that actually $\l_F =(0, +\infty )$: Let $x_0 >0 $, let G be such that for any $x>0 $:
$$ G \left( \Xpos_{\cdot \wedge \spos_x }\right) = \un_{(x_0 , \infty )} (x)  F\left( \Xpos_{\cdot \wedge \spos_{x_0}} \right)\un_{ \{ 
\Xpos_{\spos_{x_0} } = x_0 \} } \; .$$
Observe that on $\{ \Xpos_{\spos_{x_0} } = x_0 \} $, we have for any $x>x_0$
$$\spos_x = \sigma_{x-x_0} ( \Xpos \circ \theta_{\spos_{x_0}} ) + \spos_{x_0} \; .$$ 
Then, Theorem \ref{williamsaut} $(ii)$ (already proved) implies that
$$  \E \left[ G \left( \Xpos_{\cdot \wedge \spos_x }\right) ; \Xpos_{\spos_x } =x \right] =
d^* u^* (x-x_0 ) \E \left[  F\left( \Xpos_{\cdot \wedge \spos_{x_0} } \right) ; \Xpos_{\spos_{x_0} } = x_0 \right] \; .$$
Since $\l_G $ is a set of full Lebesgue measure we can assume that $x>x_0$ is in $\l_G $ and consequently
$$  \E \left[ G \left( \Xpos_{\cdot \wedge \spos_x }\right) ; \Xpos_{\spos_x } =x \right]=
\E \left[ G \left( \widehat{X}^{\tau_x } \right) ; X_{\tau_x } =x \right] \; .$$
But 
$$ \un_{\{  X_{\tau_x } =x \} }  G \left( \widehat{X}^{\tau_x } \right) =  
F \left( \widehat{Z}^{\tau_{x_0} (Z)}  \right) \un_{ \{ Z_{\tau_{x_0} (Z) }=x_0 \;  {\rm and} \;  X_{\tau_{x-x_0}} = x-x_0 \} } \; .$$
where $Z= X\circ \theta_{\tau_{x-x_0 } } $. Applying the Markov property at time $\tau_{x-x_0 } $, we get that 
$$  \E \left[ G \left( \Xpos_{\cdot \wedge \spos_x }\right) ; \Xpos_{\spos_x } =x \right]= d^* u^* (x-x_0 ) 
\E \left[ F \left( \widehat{X}^{\tau_{x_0}} \right) ; X_{\tau_{x_0}}=x_0 \right] $$
which implies the desired result. \cqfd

\subsection{Bismut's decomposition.}

As a consequence of Theorem \ref{williamcont}
and Lemma \ref{pseudobismut}, we extend to real L\'evy processes Bismut's decomposition
of the excursion above the infimum.
\begin{theorem}
\label{truebismut}
Assume that $ d^*$ is positive.
Then for any nonnegative measurable functionals $G$ and $D$ on $\OO $
 and any nonnegative measurable function $f$,
\begin{multline*}
N \left( \int_0^{\zeta } ds \; G\left( \o_{\cdot \wedge s } \right) f(\o_s )
D\left( \o \circ \theta_s \right) \right) \; = \\
\int_0^{+\infty } dx \; f(x)u^*(x) \EE \left[ G\left( \Xpos_{\cdot \wedge \spos_x }
\right) \arrowvert \Xpos_{\spos_x } =x \; \right] 
\EE \left[ D\left( X_{\cdot \wedge \tau_{-x } } \right) \right] \; .
\end{multline*}
\end{theorem}

\vspace{7mm}

\rem The spectrally positive case is due to Chaumont (see 
\cite{Ch2} ). 

\vspace{7mm}

\dem
Apply Markov property under $N$ in order to get:
$$ N \left( \int_0^{\zeta } ds \; G\left( \o_{\cdot \wedge s } \right) f(\o_s )
D\left( \o \circ \theta_s \right) \right) \; =\; 
N \left( \int_0^{\zeta } ds \; G\left( \o_{\cdot \wedge s } \right) f(\o_s )
d( \o_s ) \right) \; , $$
\noi where, for any positive number $x$,  
$ d(x)$ stands for $ \EE \left[ D\left( X_{\cdot \wedge \tau_{-x } } \right) \right] \; $.
Then, by Proposition \ref{bismut}, we have
$$ N \left( \int_0^{\zeta } ds \; G\left( \o_{\cdot \wedge s } \right) f(\o_s )
D\left( \o \circ \theta_s \right) \right) \; =\;
\EE \left[ \int_0^{L^*_{\infty } } du \; G\left( \wh X^{\locc_u } \right)
f( X_{\locc_u } ) d( X_{\locc_u } ) \right] $$ 
\noi and we use Lemma \ref{pseudobismut} and Theorem \ref{williamcont}
to complete the proof.
\cqfd

\vspace{7mm}

We have seen in Section 3 that the excursion under the supremum 
may end with a jump if  $\pi $ charges $(0, +\infty )$. Theorem \ref{revexc} provides a reversion formula for the excursion 
under the supremum  
at its final jump time. If we assume that $d^* $ is positive, then, 
the excursion may end continuously, as to say $ \o_{\zeta } =0 $. 
More precisely it is clear that $N^* (\o_{\zeta} =0 ) =0 $ if $d^*=0 $; let us show that 
$N^* (\o_{\zeta} =0 )=+\infty $ if $d^* >0 $:
\ba
N^* \left( 1-e^{-\l \zeta } \; ;\; \o_{\zeta }= 0 \right) & = & \l \int_0^{+\infty } ds e^{-\l s } N^* \left( \zeta > s\; ;\; \o_{\zeta } =0 \; \right) \\
& = & \l \int_0^{+\infty } ds e^{-\l s } N^* \left( [ \PP (X_{\tau_{a} } =a )]_{a=-\omega_s }\; ;\; \zeta > s\;  \right) \\
& = & \l \int_0^{+\infty } ds e^{-\l s } N^* \left( d^* u^* ( -\o_s  ) \; ;\; \zeta > s \; \right) \; .
\ea
By a change of variable, we have 
$$ \l \int_0^{+\infty } ds e^{-\l s } N^* \left( d^* u^* ( -\o_s  ) \; ;\; \zeta > s \; \right) =
\int_0^{+\infty } ds e^{-s } N^* \left( d^* u^* ( -\o_{s/ \l }  ) \; ;\; \zeta > s/ \l  \; \right) $$
Then for any $\l >0 $, we have
$$ N^* \left( 1-e^{-\l \zeta } \; ;\; \o_{\zeta }= 0 \right) \geq \frac{d^*}{e} 
\inf_{x\in (0, 1] } u^* (x) \;   N^* \left( \sup_{s\in [0, 1/\l ]} (-\o_s ) 
\leq 1 ;\; \zeta > 1/ \l \right) \; .$$ 
But 
\begin{multline*}
\lim_{\l \rightarrow +\infty }  N^* \left( 1-e^{-\l \zeta } \; ;\; \o_{\zeta }= 0 \right) =  
N^* \left( \o_{\zeta }= 0 \right) \\ \geq \frac{d^*}{e} \inf_{x \in (0, 1] } u^* (x) \;   
\lim_{\l \rightarrow +\infty } N^* \left( \sup_{ s\in [0, 1/\l ]} (-\o_s  ) 
\leq 1 ;\; \zeta > 1/ \l \right) = +\infty \; . 
\end{multline*}

\vspace{7mm}

The following theorem complements Theorem \ref{revexc} by providing 
a reversion identity for the excursion 
under the supremum ending continuously.

\begin{theorem}
\label{revexccont}
Assume that $d d^* >0 $. Then, for any nonnegative measurable functional $F $ on $\OO $,
$$ d\; N^* \left( F(\wh \o^{\zeta } ) \; ;\; \o_{\zeta } = 0 \; \right)\; =\; d^* \;  N\left( F(\o ) \; ;\; \o_{\zeta } = 0 \; \right) \; . $$
\end{theorem}

\vspace{7mm}
  
\rem The theorem remains true if $d^* d =0 $: in that case, it just means that 
either $N^* (\o_{\zeta } =0 )=0 $ or $N( \o_{\zeta } =0 )=0  $.

\vspace{7mm}

\dem We prove the following identity:
\begin{equation} 
\label{decompbismut1}
d\; N^* \left( \int_0^{ \zeta } ds F ( \wh \o^{\zeta }_{\cdot \wedge s} ) \; ;\; \o_{\zeta } = 0 \; \right) \; =\; 
d^* \; N \left( \int_0^{\zeta } ds F( \o_{\cdot \wedge s} ) \; ;\; \o_{\zeta } =0 \; \right) \; , 
\end{equation}
\noi which easily leads to the statement of the theorem. First, observe that 
$$ \widehat{\o }^{\zeta }_{\pw s }= \widehat{\o \circ \theta_{\zeta -s }} \; .$$
After the change of variable $s\rightarrow \zeta -s $, the Markov property under $N^* $ combined with the latter observation 
give
$$ d\; N^* \left( \int_0^{ \zeta } ds F ( \wh \o^{\zeta }_{\cdot \wedge s} ) \; ;\; \o_{\zeta } = 0 \; \right) = 
d\; N^* \left( \int_0^{ \zeta } ds \EE \left[  F\left( \widehat{X}^{\tau_{-\o_s } } \right) \;  ;\; 
X_{\tau_{-\o_s } }= -\o_s \right] \right) \; .$$
By Proposition \ref{bismut} and the dual version of Lemma \ref{pseudobismut}, it follows that 
\ba
d\; N^* \left( \int_0^{ \zeta } ds F ( \wh \o^{\zeta }_{\cdot \wedge s} ) \; ;\; \o_{\zeta } = 0 \; \right) &=& 
d \int_0^{+\infty } dx \; u (-x) 
\EE \left[ F\left( \widehat{X}^{\tau_{x} } \right) \;  ;\; 
X_{\tau_{x} }= x \right] \\
&=& dd^*  \int_0^{+\infty } dx \; u^* (x)  u (-x) 
\EE \left[ F\left( \widehat{X}^{\tau_{x} } \right) \mid 
X_{\tau_{x} }= x \right] \; .
\ea 
Use Proposition \ref{bismut} and Lemma \ref{pseudobismut} to get 
$$ dd^*  \int_0^{+\infty } dx \; u^* (x)  u (-x) 
\EE \left[ F\left( \widehat{X}^{\tau_{x} } \right) \mid 
X_{\tau_{x} }= x \right] = d^* \; N \left( \int_0^{\zeta } ds F( \o_{\cdot \wedge s} ) \; ;\; \o_{\zeta } =0 \; \right) $$
that is the desired result. \cqfd 

\subsection{The second Williams' decomposition theorem.}

Williams has shown in \cite{Wil74} that the 
Brownian excursion splits at its maximum in two three-dimensional Bessel processes
stopped at a certain hitting time. In this section, we extend this result to 
general L\'evy processes. To simplify, we set
$$ Z^T = \XdT \circ \theta_{\tdT } \qquad {\rm under} \qquad \PP (\cdot \cap  \{ \tdT <\infty \} ) \qquad {\rm and } \qquad 
Z= \Xpos \circ \theta_{\tpos_x  } \; .$$
We first prove the following proposition.
\begin{proposition}
\label{ind}

\noi (i) For any bounded measurable functional $F$ on $\OO $, 
$$ \EE \left[ F\left( \Xpos_{\cdot \wedge \tpos_x } 
\right) \right] \; =\; N\left(  F \left( \o_{\cdot \wedge \tau_x (\o ) }
\right) \U \left( ( -\o_{\tau_x (\o ) } , 0 ] \right) ;\; \tau_x (\o ) <\infty \right) \; . $$
Thus, $\PP ( \Xpos_{\tpos_x } =x ) >0 $ if and only if $d^*> 0$.

\noi (ii) If $d^*>0 $, then, under $\PP (\cdot \mid  \Xpos_{\tpos_x } =x ) $, 
$\underleftarrow{Z} $ and $\underrightarrow{Z}$ are mutually independent; the process $\Zd$ 
has the same law as $\Xpos $ and the law of $\Zg $ is characterized 
by the following identity that holds for any nonnegative measurable functional 
$F$ on $\OO $:
$$ \EE \left[ F( \Zg )\right] = \frac{1}{\u ((-x , 0]) }\EE \left[ 
\int_0^{L_{\infty }}dv F( X_{\pw L^{-1}_v } ) \un_{ \{ X_{L^{-1}_v } >-x 
\} } \right] \; .$$

\end{proposition}

\vspace{7mm}

\rem Recall from 
Section 2.1, that under $\PP (\cdot \mid \Xpos_{\tpos_x  } =x )$, the process $x+Z$ 
is markovian with a transition kernel given by
$$ p^+_t (y, dz ) = \frac{\u ( (-z, 0] )}{ \u ((-y, 0] ) } q^+_t (y, dz ) 
\; , \quad y\geq 0 \; , $$
where $q^+_t $ stands for the semigroup of the L\'evy process killed in 
$(-\infty , 0] $. Then, the latter proposition combined with 
Lemma \ref{pseudobismut} gives the following corollary. 
\begin{corollary}
\label{markov}
Assume that $dd^* >0 $. Let $\Xpos (x) $ denote the L\'evy process started at 
$x>0 $ and conditioned to stay positive. The path $\Xpos (x)$ has 
the following 
decomposition at its infimum:

\noi (i) The pre-infimum process $\underleftarrow{\Xpos} (x) $ and 
the post-infimum process $\underrightarrow{\Xpos} (x) $ are mutually independent
and $\underrightarrow{\Xpos} (x) $ is distributed as the L\'evy process conditioned to stay positive started at $0$.

\noi (ii) The law of the infimum of $\Xpos (x)$ admits a density with 
respect to Lebesgue measure that is given by: 
$$ y\; \longrightarrow \un_{ [0,x]} (y) \frac{ u(y-x)}{\u ((-x, 0] )} $$
and under $\PP ( \cdot \mid \inf \Xpos (x) =y ) $, 
$\underleftarrow{\Xpos} (x) $ is distributed as $x + X_{\pw \tau_{y-x} } $
under $\PP ( \cdot \mid X_{\tau_{y-x }} = y-x )$. 

\end{corollary}

\vspace{7mm}

\rems (i) We can actually show that the corollary remains true even if $ d^* =0 $.

\noi (ii) When $X$ does not drift to $-\infty $ the result is Theorem 5 of Chaumont \cite{Ch2}
(see also \cite{Ch} and see also \cite{Ch1}). 
   
\vspace{7mm}

\noi {\bf Proof of Proposition \ref{ind}.}
Let $T$ be independent of $X$ and exponentially distributed with parameter 
$\alpha $. Recall from Lemma \ref{independance} that for any nonnegative 
measurable functional $H$ defined on $\OO $, 
$$ \EE \left[  H\left( \XdT \right) \right] \; =\; \frac{\alpha}{\kappa (\alpha , 0) } 
N\left( \int_0^{\zeta } ds \; e^{-\alpha s } H\left( \o_{\cdot \wedge s } \right)
\right) \; . $$
Let $G$ and $D$ be nonnegative measurable functionals on $\OO $. Take 
$$ 
H(\o ) = \un_{\{ \tau_x (\o ) < \infty \} } 
G( \o (\pw \tau_x (\o ) )) D( \o \circ \theta_{\tau_x (\o ) } ) \; , 
\quad \o \in \OO \; .$$
Then,
\begin{multline*}
\EE \left[ G \left( \XdT_{\cdot \wedge \tdT } \right) D\left( \XdT \circ \theta_{\tdT } \right) ;\; \tdT <\infty \right]  \\
=\frac{\alpha}{\kappa (\alpha , 0) } N\left(  \un_{\{ \tau_x (\o ) < \infty \} } 
G\left( \o_{\cdot \wedge \tau_x (\o ) } \right) e^{-\alpha \tau_x (\o ) } \int_0^{\zeta -\tau_x (\o )} ds \; e^{-\alpha s }
D\left( \left( \o \circ \theta_{\tau_x (\o ) } \right)_{\cdot \wedge s} \right) \right) . 
\end{multline*}
\noi Apply the Markov property under $N$ in order to get
\begin{multline}
\label{ind1}
\EE \left[ G \left( \XdT_{\cdot \wedge \tdT } \right) D\left( \XdT 
\circ \theta_{\tdT } \right) ; \;  \tdT <\infty \right] 
\\ =  
N\left(  \un_{\{ \tau_x (\o ) < \infty \} } e^{-\alpha \tau_x (\o ) } G\left( \o_{\pw \tau_x (\o ) } \right)
d_{\alpha } ( \o_{\tau_x (\o )} ) \right)   \; ,
\end{multline}
where for any $a>0 $, 
$$ d_{\alpha } (a)= \frac{1}{\kappa (\alpha , 0) }
\EE \left[ D( X_{ \pw T} ) 
\; ;\; I_T > -a \right]\; .$$
Take $D=1 $ in (\ref{ind1}). We get  
$$ \EE \left[ G \left( \XdT_{\cdot \wedge \tdT } \right) \right]
= N\left(  \un_{\{ \tau_x (\o ) < \infty \} }
e^{-\alpha \tau_x (\o ) } G\left( \o_{\pw \tau_x (\o ) } 
\right) \frac{\PP (I_T > -\o_{\tau_x (\o ) } )}{\kappa (\alpha , 0) } \right) \; .$$
Next we need the following lemma:
\begin{lemma}
\label{lien}
Let $\alpha >0 $ and $T_{\alpha }$ be 
independent of 
$X$ and exponentially distributed with parameter $\alpha $. Let $F$, $G$ and 
$K$ be three bounded measurable functionals on $\OO $.
Under the same assumptions as Theorem \ref{williamsaut}, we have 
$$ \lim_{\alpha \rightarrow 0} \EE \left[ 
F\left( \Xd^{T_{\alpha }}_{\cdot \wedge \sd^{T_{\alpha } }_x } \right)
G( \Yg^{ T_{\alpha }} ) K( \Yd^{T_{\alpha } } ) ; \;  \sd^{T_{\alpha } }_x <\infty \right]=
\EE \left[ F\left( \Xpos_{\pw \spos_x } \right) G( \Yg ) K(\Yd ) \right]
\; . $$
\end{lemma}
\dem It is sufficient to show the limit
$$ \lim_{s\rightarrow +\infty } \EE \left[ 
F\left( \Xd^{s}_{\cdot \wedge \sd^{s }_x } \right)
G( \Yg^{ s} ) K( \Yd^{s } ) ;\; \sd^{s }_x <\infty \right]=
\EE \left[ F\left( \Xpos_{\pw \spos_x } \right) G( \Yg ) K(\Yd ) \right]
\; . $$
We use the notation of Section 2.2. Recall that 
$ A^+_s = \int_0^s du \un_{\{ X_u >0 \} } $. From Theorem
\ref{bertoinresult}, we have 
\begin{equation}
\label{lien1} 
(\Xd^s_u)_{0 \leq u <s-\underline{g}_s } \; \overset{(law)}{=} \;
(\Xpos_u)_{0\leq u < A^+_s } \; . 
\end{equation}
Set $\beta = \spos_x + \underline{g} (Y) $. By (\ref{limite}), 
$\beta < \infty $. In order to avoid cumbersome notations, we denote by 
$W$ and $W'$ respectively the pre-infimum process and the post-infimum 
process of $\Xpos \circ \theta_{ \sigma_x ( \Xpos , A^+_s ) }$. 
Let $M>0 $ be an upper bound for $F, G $ and $K$. Observe that on $\{ A^+_s > \beta \}$, we have 
$ \sigma_x ( \Xpos , A^+_s ) =\spos_x $ and 
\begin{equation}
\label{lien2}
F \left( \Xpos_{\cdot \wedge \sigma_x ( \Xpos , A^+_s ) } \right) 
G\left( W\right)
K\left( W' \right) =
F\left( \Xpos_{\pw \spos_x } \right) G( \Yg ) K(\Yd )\; .
\end{equation}
Then by (\ref{lien1}) and (\ref{lien2}), we have
\begin{multline*}
\EE \left[ 
 F\left( \Xd^{s}_{\cdot \wedge \sd^{s }_x } \right)
G( \Yg^{ s} ) K( \Yd^{s } ) \right] = \EE \left[ F \left( \Xpos_{\cdot \wedge \sigma_x ( \Xpos , A^+_s ) } \right) 
G\left( W \right) K\left( W' \right) ; \;  A^+_s \leq  \beta \right] \\
+ \EE \left[ F\left( \Xpos_{\pw \spos_x } \right) G( \Yg ) K(\Yd ) ; \;  A^+_s > \beta \right] \; .
\end{multline*}
Consequently,
$$ \left| \;  \EE \left[ 
 F\left( \Xd^{s}_{\cdot \wedge \sd^{s }_x } \right)
G( \Yg^{ s} ) K( \Yd^{s } ) \right] - 
 \EE \left[ F\left( \Xpos_{\pw \spos_x } \right) G( \Yg ) K(\Yd ) \right] \; 
\right| \leq \; 2 M^3 \; \PP ( \beta \geq A^+_s ) \; . $$
But $\; \lim_{s\rightarrow \infty } A^+_s = \infty $. So 
$\lim_{s\rightarrow \infty } 
\PP ( \beta \geq A^+_s ) =0 $ which yields the lemma. \cqfd

\vspace{5mm}

Let us achieve the proof of the proposition: Lemma \ref{lien} implies  
$$ \lim_{\alpha \rightarrow 0} 
\EE \left[ G \left( \XdT_{\cdot \wedge \tdT } \right) \right] =
\EE \left[ G \left( \Xpos_{\cdot \wedge \tpos_x  } \right) \right] \; $$
Next, we deduce from Lemma \ref{independance} that for any $a>0 $
$$ \lim_{\alpha \rightarrow 0 } \frac{\PP (I_T > - a )}{\kappa (\alpha , 0) }
= \EE \left[ \int_0^{L^{\infty }} dv \un_{ \{ X_{L^{-1 }_v} > -a \} } \right] =
\u ( (-a , 0] ) \; ,$$
which yields $(i)$ by dominated convergence.

  Assume now that the drift coefficient $d^*$ is positive. Let $D_1 $ and $D_2 $ be 
two nonnegative measurable functionals on $\OO$. Take, 
$$ D( \o ) = D_1 ( \underleftarrow{\o } ) D_2 ( 
\underrightarrow{\o }) \qquad {\rm and} \qquad   G(\o ) = \un_{ \{ \tau_x (\o ) <\infty ;\; 
 \o (\tau_x (\o )) = x \} } \; ,$$
in (\ref{ind1}). From Lemma \ref{independance}, we note that 
\begin{eqnarray*}
d_{\alpha } (a) & = & \frac{1}{\kappa (\alpha , 0) }
\EE \left[ D_1 ( \XgT ) \; ; \; I_T >-a \right] \EE 
\left[ D_2 ( \XdT ) \right]  \\
&=&  \EE \left[ \int_0^{L_{\infty } } dv e^{-\alpha L^{-1 }_v }
D_1 ( X_{\pw L^{-1 }_v } ) \un_{ \{ X_{L^{-1}_v } >-a \} } \right]
\EE \left[ D_2 ( \XdT ) \right] \; .
\end{eqnarray*} 
Thus, (\ref{ind1}) gives 
\begin{multline} 
\label{ind2}
\EE \left[ F \left( \XdT_{\cdot \wedge \tdT } \right) 
D\left( \XdT \circ \theta_{\tdT } \right) ;  \; \tdT <\infty 
;\; \XdT_{\tdT }=x \right] =\EE \left[ D_2 ( \XdT ) \right]\\ 
\times  N\left(e^{-\alpha \tau_x (\o ) } \; ;\;  \tau_x (\o )<\infty ;\; 
\o_{\tau_x (\o ) } =x \right) 
 \EE \left[ \int_0^{L_{\infty } } dv e^{-\alpha L^{-1 }_v }
D_1 ( X_{\pw L^{-1 }_v } ) \un_{ \{ X_{L^{-1}_v } >-x \} } \right]
 \; .
\end{multline}
To get $(ii)$, pass to the limit $\alpha \rightarrow 0 $ in (\ref{ind2}) using 
Lemma \ref{lien} to write:
\begin{multline*}
\lim_{\alpha \rightarrow 0} 
\EE \left[D\left( \XdT \circ \theta_{\tdT } \right) ;\; \tdT <\infty  ;\; \XdT_{\tdT }=x \right] 
\; = \\
\EE \left[ F( \Xpos_{\pw \tpos_x  } ) D_1 ( \Zg ) D_2 (\Zd ) \; ;\; 
\Xpos_{\tpos_x  } = x \right] \; .
\end{multline*}
and 
$$ \lim_{\alpha \rightarrow 0} \; \EE \left[ D_2 ( \XdT ) \right]
\; = \; \EE \left[ D_2 (\Xpos ) \right] \; .$$
\cqfd

\vspace{7mm}

We are now able to state the second Williams' decomposition theorem.
\begin{theorem}
\label{williamdeux}
Assume that $\pi $ charges $(0, +\infty ) $. Suppose also that $d^* >0$ 
and that $X$ oscillates.

\noi (i) The law of $\o_{\ggg (\o ) } $ under $N$ admits a density with 
respect to Lebesgue measure that is given by 
$$ x \longrightarrow \frac{1}{d^*} N^* \left( 
\tau_{-x} (\o ) <\infty \right) N\left( \tau_x (\o ) < \infty  ; \;  \o_{\tau_x (\o ) } =x \right) \; .$$

\noi (ii) Under $N(\cdot \mid \o_{\ggg (\o ) } =x )$, the processes 
$ \o_{ \pw \ggg (\o ) } $ and $ \o \circ \theta_{\ggg (\o ) } $ are mutually
independent. Furthermore,

- the process $\o_{\pw \ggg (\o ) } $ is distributed as $\Xpos_{\pw \tpos_x  } $ under
$\PP (\cdot \mid \Xpos_{\tpos_x  } =x ) $;

- the law of $ \o \circ \theta_{\ggg (\o ) } $ is absolutely continuous 
with respect to the law of $\Xneg_{\pw \tneg_{-x } } $ and the corresponding density 
is $\varphi (\Xneg_{\tneg_{-x }} ) $, where 
$$ \frac{1}{\varphi (y)} = \u^* ([0, -y )) N^* (\tau_{-x } (\o ) <\infty ) \; , \quad y \in (-\infty , 0) \; . $$
\end{theorem}

\vspace{7mm}

\rem The spectrally positive case is due to Chaumont in \cite{Ch2} or \cite{Ch}. 

\vspace{7mm} 

\dem
Let $G$ and $D$ be two nonnegative measurable functionals on $\OO $. From 
Proposition \ref{ind} $(i)$ and the corresponding dual equality, we get for any $x>0 $,
$$\EE \left[ G\left( \Xpos_{\cdot \wedge \tpos_x } \right) \right] \; =\; N\left( G\left( \o_{\cdot \wedge \tau_x (\o ) } \right) 
\U \left( (-\o_{\tau_x (\o) } , 0] \right)
;\; \tau_x (\o) <\infty \right) $$
and 
$$ \EE \left[ D \left( \Xneg_{\cdot \wedge \tneg_{-x }} \right) \right] \; =\; N^* \left( D\left( \o_{\cdot \wedge \tau_{-x} (\o ) } \right) 
\U^* \left( [0, -\o_{\tau_{-x } (\o ) } ) \right) ;\; \tau_{-x} (\o) <\infty \right) \; .$$
The various assertions of the theorem then follow from the identity
\begin{multline}
\label{approx}
N \left( G\left( \o_{\cdot \wedge \gs (\o ) } \right) D\left( \o \circ \theta_{\gs (\o ) } \right) \right) \\
= \fr{1}{d^* }
\int_0^{+\infty } dx \;  N\left( G\left( \o_{\cdot \wedge \tau_x (\o ) } \right) \; ;\; \o_{\tau_x (\o) }=x \right)
N^* \left( D\left( \o_{\cdot \wedge \tau_{-x} (\o ) } \right) ;\; \tau_{-x} (\o) <\infty \right) \; , 
\end{multline}
which we now prove.

  Let $a$ and $b$ be two positive real numbers. Observe that
$$ \{ a\leq \o_{\gs (\o ) } <a+b \} \; =\; \{ \tau_a (\o ) <\infty \} \cap  
\{ \sup \; \o \circ \theta_{\tau_a (\o ) } \; < b+a -\o_{\tau_a (\o ) } \; \} \; . $$
\noi On this event, $ \o \circ \theta_{\gs (\o ) } $ is 
the post-supremum process of $\o \circ \theta_{\tau_a (\o) } $. Hence, by the 
Markov property under $N$ at $\tau_a (\o) $, we have   
\begin{multline*}
 N\left( G\left( \o_{\cdot \wedge \tau_a (\o ) } \right) D\left( \o \circ \theta_{\gs (\o) } \right) \; ;\; a\leq \o_{\gs (\o) } < a+b \; \right)
\\ = N\left( G\left( \o_{\cdot \wedge \tau_a (\o ) } \right)  d( \o_{\tau_a (\o) } ) ;\; 
 \tau_a (\o ) <\infty\right) \; , 
\end{multline*}
\noi where for any positive number $x$, 
$$ d(x) \; =\;  \EE \left[ D\left( \Xupd^{\tau_{-x } } \right) \; ;\;  S_{\tau_{-x}} < b+a -x \; \right] \; . $$
\noi Let us write $d(x)$ in a more suitable form. First, observe that 
$$ d(x) \; =\; \EE \left[ \sum_{j\in \i^* } \un_{\{ X_{g_j }< b+a-x \; ;\; 
I_{g_j } >-x \} }
D \left( \o^j_{ \cdot \wedge \tau_{-x-X_{g_j }} ( \o^j ) } \right) \un_{\{ \tau_{-x-X_{g_j }} ( \o^j ) <\infty \} }
\right] 
\; . $$
\noi Then, apply the compensation formula to get
$$ d(x) = \EE \left[ \int_0^{\infty } du 
\un_{\{ X_{\locc_u } <b+a-x \; ;\; I_{\locc_u } >-x \} }
N^* \left( D\left( \o_{\cdot \wedge \tau_{-x-X_{\locc_u }} (\o ) } 
\right) \un_{ \{ \tau_{-x-X_{\locc_u }} (\o ) <\infty \} } 
\right) \right] .$$
\noi Set for any real numbers $x$ and $y$, 
$$ v(x,y)\; =\; 
\un_{ [0, +\infty ) } (x) \un_{ [0, +\infty ) } (y)  u^* (y ) \PP \left( I_{\tau_y } >-x \; |\; X_{\tau_y } =y \; \right) \; . $$
\noi Then, for any $a\leq x <a+b $, Lemma \ref{pseudobismut} implies
\ba
d(x) &=& \int_0^{a+b-x} dy u^* (y) \EE \left[ \un_{\{ I_{\tau_y } > -x \} }  
N^* \left( D ( \o_{\cdot \wedge \tau_{-x-y} (\o ) } ) ;\; \tau_{-x-y} (\o ) <\infty \right) \mid X_{\tau_y } = y \right] \\
&=&\; \int_a^{a+b } dy \; 
v(x, y-x ) N^* \left( D ( \o_{\cdot \wedge \tau_{-y} (\o ) } ) ;\; \tau_{-y} (\o ) <\infty \right) \; .
\ea
\noi Thus, 
\begin{multline}
\label{approxox}
 N\left( G\left( \o_{\cdot \wedge \tau_a (\o ) } \right) D\left( \o \circ \theta_{\gs (\o) } \right) \; ;\; a\leq \o_{\gs (\o) } 
< a+b \; \right) 
=\int_a^{a+b } dy  \\
\times N \left( \un_{ \{ \tau_a (\o) <\infty \} }
v(\o_{\tau_a (\o) }, y-\o_{\tau_a (\o ) } )  G\left( \o_{\cdot \wedge \tau_a (\o ) } \right) \; \right) N^* 
\left( D\left( \o_{\cdot \wedge \tau_{-y} (\o ) } ;\; \tau_{-y} (\o ) <\infty \right) \right) 
\end{multline}
Next, set for any positive integer $n$,   
$$ \; m_n = \frac{[2^n \o_{\gs (\o ) } ]}{2^n} \qquad {\rm and }\qquad 
y_n =\frac{[2^n y]}{2^n}\; , \quad y\geq 0 
\; .$$
We apply (\ref{approxox}) with $a=i2^{-n } $ and $ b= 2^{-n } $ for every integer $i\geq 0$, and we sum over $i$. It follows that  
\begin{multline}
\label{aproprox}
N\left( G\left( \o_{\cdot \wedge \tau_{m_n } (\o ) } \right) D \left( \o \circ \theta_{\gs (\o )} \right) \right)\; = \;
\int_0^{+\infty } dy \\ 
\times N \left( \un_{ \{ \tau_{y_n } (\o) <\infty \} }
v(\o_{\tau_{y_n } (\o ) } , y - \o_{\tau_{y_n } (\o ) } )  G\left( \o_{\cdot \wedge \tau_{y_n } (\o ) } \right)  \right)
N^* \left( D\left( \o_{\cdot \wedge \tau_{-y} (\o ) } ;\; \tau_{-y} (\o ) <\infty \right) \right)  . 
\end{multline}

  Let $\eps , A >0 $. It is sufficient to prove (\ref{approx}) for  
$$ G(\o ) = \un_{ \{ \zeta (\o ) > \eps ;\; \sup \o \leq  A \}  } 
F( \o ) \qquad {\rm and} \qquad  D(\o ) = \un_{ \{ \zeta (\o ) > \eps  \} } K( \o ) $$
where
$F(\o ) $ and $K(\o )$ depend continuously on the values of $\o $ 
at some finitely many positive times. 
Assumption {\bf (A)} implies that $X$ attains continuously its supremum on any finite time interval
(see Millar \cite{Mil77}). So does 
the excursion above the infimum. 
Thus, $\tau_{m_n } (\o ) $ increases  to $\gs (\o) $
when $n$ goes to infinity $N$-almost everywhere. Then, 
$$ \lim_{n\rightarrow \infty }  \un_{ \{ \tau_{m_n } (\o) > \eps ;\; m_n \leq A \} } =  \un_{ \{ \gs (\o)  > \eps ; \; 
\sup \o \leq A  \} } \; , \quad N-{\rm a.e.} \; .$$
Since $N (  \gs (\o)  > \eps ;\; \sup  \o \leq A  ) \leq N (  \zeta (\o)  > \eps ) <\infty $, 
dominated convergence applies in the left side of (\ref{aproprox}) and we get:  
\begin{multline}
\label{approx3}
\lim_{n\rightarrow \infty } N\left( G\left( \o_{\cdot \wedge \tau_{m_n } (\o ) } \right) D \left( \o \circ \theta_{\gs (\o )} \right) 
\right) \\
=N\left( F\left( \o_{\cdot \wedge \gs (\o ) } \right) 
K \left( \o \circ \theta_{\gs (\o )} \right) ;\; \sup \o \leq A  ;\; (\zeta -\gs (\o ) )\wedge  \gs (\o ) > \eps \right)  \; .
\end{multline}

 We now turn to the limit of the right hand side of (\ref{aproprox}): Recall
(\ref{rog2}) from Section 2.1 
$$ \PP \left( \exists t \in (0, +\infty ) \; : \; S_{t- } = X_{t- } < X_t \; \right) \;= \; 0 \;  . $$
\noi It implies that for any positive number $y\; $, $\; N\left( \o_{\tau_y (\o )- } =y < \o_{\tau_y (\o ) } \; \right) \; =\; 0 \; $. 
We also recall that $\o (0) = 0 $, $N $-a.e. . Thus,
\begin{equation}
\label{precision1}
N-{\rm a.e.} \qquad \lim_{n\rightarrow \infty } 
\un_{ \{ \o_{ \tau_{y_{n}} (\o )} \leq y \} } = \un_{ \{ \o_{\tau_y (\o ) } = y 
\} } \; .
\end{equation}
Next, for any $x>0 $, 
$ \lim_{\eps \rightarrow 0 } v(x, \eps ) = u^*(0+)= \frac{1}{d^* }$ because 
\begin{equation*}
0\leq d^* u^*( \eps ) - d^* v(x, \eps )= d^* u^*(\eps ) \PP \left( I_{\tau_{\eps } } \leq -x \mid X_{\tau_{\eps } }=\eps \right)
\leq \PP \left( I_{ \tau_{\eps } } \leq  -x \right) \xrightarrow[\eps \rightarrow 0 ]{\quad } 0 \; .  
\end{equation*}
Thanks to (\ref{precision1}) we get $N$-a.e. 
$$ \lim_{n\rightarrow \infty } \un_{\{ \tau_{y_n } (\o ) <\infty \} }
v\left( \o_{\tau_{y_n } (\o ) } , y -\o_{\tau_{y_n } (\o ) } \right) \; =\; u^* (0+) \un_{\{ \o_{\tau_y (\o ) } =y \} } =
\frac{1}{d^* } \un_{\{ \o_{\tau_y (\o ) } =y \} } \; . $$
\noi Observe that 
$N$-a.e.  on $\{ \o_{\tau_y (\o ) } = y \} $, $ \tau_{y_n } (\o ) $ increases towards 
$ \tau_y (\o ) $. Therefore, by dominated convergence 
\begin{multline*}
 \lim_{n\rightarrow \infty } N\left(  \un_{\{ \tau_{y_n } (\o ) <\infty \} }
v(\o_{\tau_{y_n } (\o ) } , y - \o_{\tau_{y_n } (\o ) } )  
G \left( \o_{\cdot \wedge \tau_{y_n } (\o ) } \right)  \right) \\
= \fr{1}{d^* } \un_{(0, A] } (y) N\left( F\left( \o_{\cdot \wedge \tau_y (\o ) } \right) ;\; \eps < \tau_y (\o ) < \infty  
;\; \o_{\tau_y (\o ) } =y \; \right) \; . 
\end{multline*}
However, for any $y>0 $, we have 
\begin{multline*}
N \left( \un_{ \{ \tau_{y_n } (\o) <\infty \} }
v(\o_{\tau_{y_n } (\o ) } , y - \o_{\tau_{y_n } (\o ) } )  G\left( \o_{\cdot \wedge \tau_{y_n } (\o ) } \right)  \right)
N^* \left( D\left( \o_{\cdot \wedge \tau_{-y} (\o ) } ;\; \tau_{-y} (\o ) <\infty \right) \right) \\
\leq  M^2 \un_{ (0, A] } (y) \sup_{x\in (0, 2^{-n }]} u^* (x)  \; N(\zeta > \eps ) N^*(\zeta > \eps ) \; ,
\end{multline*}
where $M$ is a bounding constant of $F$ and $K$. Then, dominated convergence applies in the right side of (\ref{aproprox}), 
which yields the desired identity  (\ref{approx}) thanks to (\ref{approxox}). \cqfd

\bibliographystyle{plain}

\end{document}